\documentclass[12pt]{amsart}
\usepackage{amssymb}
\usepackage{hyperref}
\usepackage{amsmath}
\usepackage{amscd}
\usepackage{graphicx}

\DeclareMathSymbol{\shortminus}{\mathbin}{AMSa}{"39}

\newcommand{\group}[1]{\mathrm{#1}}

\newcommand{\rep}[1]{\rep{#1}}

\newcommand{\Tr}{\operatorname{Tr}}

\newcommand{\End}{\operatorname{End}}

\newcommand{\Cat}{\operatorname{Cat}}

\newcommand{\Cay}{\operatorname{Cay}}

\newcommand{\C}{\mathbb{C}}

\newcommand{\N}{\mathbb{N}}

\newcommand{\R}{\mathbb{R}}

\newcommand{\V}{\mathbf{V}}
\newcommand{\W}{\mathbf{W}}

\newcommand{\Y}{\mathsf{Y}}

\newtheorem{thm}{Theorem}[section]

\newtheorem{cor}[thm]{Corollary}

\theoremstyle{definition}
\newtheorem{definition}[thm]{Definition}

\title[$YM_2$/Hurwtiz Correspondence]{On the 2D Yang-Mills/Hurwitz Correspondence}

\author{Jonathan Novak}

\begin{document}


\maketitle
\tableofcontents
\pagebreak

\section{Introduction}

	\subsection{Graph theory}
		Let $\group{G}$ be a connected finite simple graph with vertices $\pi,\rho,\sigma,\dots$. 
		Let $\C\group{G}$ be the free Hilbert space over the vertices of $\group{G}$
		modeled as the space of kets $|\pi \rangle,|\rho\rangle, |\sigma\rangle,\dots$.
		The adjacency operator $K \in \End \C\group{G}$ encodes the edges of $\group{G}$,
					
			\begin{equation}
				\langle\sigma| K |\rho\rangle = \begin{cases}
					1, \text{ if } \rho,\sigma \text{ adjacent}\\
					0, \text{ otherwise}.
				\end{cases}
			\end{equation}
		
		\noindent
		More generally, 
		the sequence of operators $K^1,K^2,K^3,\dots$ is meaningful in that 
		$\langle \sigma| K^r |\rho\rangle$ counts $r$-step walks $\rho \to \sigma$ in $\group{G}$.
		The exponential adjacency operator $\Psi_t = e^{-tK}$, 
		which for regular graphs is essentially the heat kernel \cite{Chung},
		encodes this information in its matrix elements,
				
			\begin{equation}
				\langle\sigma| \Psi_t |\rho\rangle = \sum_{r=0}^\infty \frac{(-t)^r}{r!}
				\langle \sigma| K^r |\rho\rangle.
			\end{equation}
			
		\noindent
		The algebraic approach to counting walks in $\group{G}$ is to find a basis of $\C\group{G}$ in which $K$ and hence $\Psi_t$ 
		act diagonally \cite{Stanley:AlCo}.
														
		Instead of asking how many walks there are between two given vertices in $\group{G}$, we may 
		ask how far apart they are. This information is encoded by the distance
		operator $D \in \End \C\group{G}$, whose matrix 
		
			\begin{equation}
				\langle \sigma | D |\rho\rangle = \mathrm{d}(\rho,\sigma)
			\end{equation}
			
		\noindent
		in the vertex basis tabulates geodesic distance $\mathrm{d}$ in $\group{G}$. 
		Spectral properties of $D$ were first explored by Graham and Lovasz \cite{GL}
		and have been much studied since; see \cite{AH}.
		
		Another way to encode distance in graphs is to use linear operators $L_1,L_2,L_3,\dots$ 
		corresponding to metric level sets,
		
			\begin{equation}
				\langle \sigma| L_r |\rho \rangle =\begin{cases}
				 1, \text{ if } \mathrm{d}(\rho,\sigma)=r \\
					0, \text{ otherwise}.
				\end{cases}
			\end{equation}

		\noindent 
		Note that $L_1=K$ and $L_r$ is zero for all $r$ larger than 
		the diameter of $\group{G}$. The exponential distance operator,
		
			\begin{equation}
				\Omega_q = \sum_{r=0}^\infty q^r L_r,
			\end{equation}
			
		\noindent
		is the entrywise exponential of $D$, 
				
			\begin{equation}
				\langle \sigma| \Omega_q |\rho \rangle = q^{\mathrm{d}(\rho,\sigma)}.
			\end{equation}
		  
		\noindent
		This matrix is natural from a statistical physics perspective:
		it tabulates the Boltzmann weight $e^{-\beta \mathrm{d}}$ corresponding to goedesic distance
		at inverse temperature $\beta =- \log q$. 
		
		The general study of exponential 
		distance operators in graph theory is just beginning; see \cite{BCLLS}.		
		To the best of this author's knowledge, the first appearance of $\Omega_q$ was in 
		Zagier's study of deformed commutation relations in quantum physics \cite{Zagier}, 
		where the exponential distance operator on the Cayley graph of the symmetric group
		as generated by the set of Coxeter transpositions $(i\ i+1)$ plays a pivotal role.
		The analysis in \cite{Zagier} culminates in an explicit formula for the determinant
		of $\Omega_q$ showing that its zeros lie on the unit circle, implying
		the existence of a Hilbert space representation of the $q$-commutation relations
		for all $-1<q<-1$, interpolating between bosons and fermions.
				
		\subsection{Hurwitz theory}
		The setting of this paper is the Cayley graph of the symmetric group 
		as generated by the full conjugacy class of transpositions, which has the 
		feature that its exponential adjacency and distance operators $\Psi_t$ and 
		$\Omega_q$ are simultaneously and explicitly diagonalizable. This is not the case 
		for the Coxeter-Cayley graph of the symmetric group.
		On the other hand, in the all-transpositions case the zeros
		of $\det \Omega_q$ are unstable: they converge on the origin as the degree of the symmetric group increases.
				
		The all-transpositions Cayley graph is intimately linked to Hurwitz theory, a classical branch of enumerative geometry
		concerned with counting maps between compact Riemann surfaces \cite{CM}, and we therefore
		refer to it as the \emph{Hurwitz-Cayley graph}. Cover counting is associated
		with the exponential adjacency operator $\Psi_t$ on this graph. For example, loop enumeration
		on the Hurwitz-Cayley graph is algebraically equivalent to computing a diagonal matrix element of 
		$\Psi_t$ and geometrically equivalent to counting simple branched covers of the Riemann sphere,
		a classical and much-studied problem \cite{DYZ}.
				
		It is a special feature of the Hurwitz-Cayley graph that
		matrix elements of the exponential distance operator $\Omega_q$ and 
		its inverse also count walks \cite{Novak:BCP}. The walks enumerated by $\langle \sigma | \Omega_q^{\pm 1}|\rho\rangle$
		are ``monotone'' with respect to an edge labeling of the Hurwitz-Cayley graph
		related to the representation theory \cite{MatNov} and order theory \cite{MakNov}
		of the symmetric group. Thus $\Omega_q$ corresponds to a monotone version of 
		Hurwitz theory \cite{GGN1,GGN2,GGN3,GGN4,GGN5} 
		which turns out to have many applications in random matrix theory and related areas.
										 		
		\subsection{Yang-Mills theory}
		The idea that two-dimensional Yang-Mills theory with $\group{U}_N$ gauge group becomes
		Hurwitz theory in the large $N$ limit \cite{Gross} has been intensively studied by physicists \cite{CMR1,CMR2,GT1,GT2,Taylor},
		who noticed that eigenvalues of the exponential adjacency operator $\Psi_t$
		of the Hurwitz-Cayley graph appear in the partition function of $YM_2$. 
		The Yang-Mills/Hurwitz correspondence works perfectly in the Calabi-Yau case, where spacetime is a torus, even manifesting
		mirror symmetry \cite{Dijkgraaf}. For non-torus spacetimes, Hurwitz theory is not an exact match for the 
		large $N$ limit of $YM_2$ due to the presence of certain ``$\Omega$-factors'' 
		which appear in the partition function with multiplicity equal to the Euler characteristic of spacetime but have no obvious 
		Hurwitz-theoretic meaning.  The large $N$ limit of $YM_2$ is necessarily a completion 
		of classical Hurwitz theory which accounts for these $\Omega$-factors. 
		
		The purpose of this paper is to point out that the $\Omega$-factors in $YM_2$ are exactly 
		the eigenvalues of the exponential adjacency operator $\Omega_q$ of the Hurwitz-Cayley graph, and that consequently
		the large $N$ limit of $YM_2$ is the union of classical and monotone
		Hurwitz theory: mixed Hurwitz theory. The mixed Hurwitz theory of a genus zero 
		target was introduced in \cite{GGN4} as a way to interpolate between results in 
		Hurwitz theory \cite{Okounkov} and random matrix theory \cite{ZJ}. It was generalized
		to higher genus targets in \cite{HIL}.
		
		%

		\subsection{Organization}
		In Section \ref{sec:HurwitzCayley} we introduce the Hurwitz-Cayley graph and
		explain what is perhaps its most important feature: a combinatorial duality 
		called the \emph{Jucys-Murphy correspondence}.
		In Section \ref{sec:YangMills} we derive the Gross-Taylor formula for the 
		chiral partition function of $YM_2$ on compact orientable surfaces \cite{CMR1,GT1},
		and explain its relationship to the graph operators $\Psi_t$ and $\Omega_q$,
		which constitutes the Yang-Mills/Hurwitz correspondence.
						
		The remaining sections of the paper illustrate the Yang-Mills/Hurwitz correspondence
		for specific choices of spacetime. This is done first 
		for the finitely many cases which can be described
		using only classical Hurwitz theory: the cylinder and 
		surfaces obtained from it by identifying boundaries or plugging holes (torus, disc, sphere). 
		We then consider the three-holed sphere and one-holed torus as representative examples 
		of the infinitely many remaining compact orientable spacetimes, where monotone Hurwitz
		theory is inescapable. In fact, the area zero limit of $YM_2$, a much-studied degeneration 
		\cite{CMR1,CMR2,Taylor}, is pure monotone Hurwitz theory. This implies a relationship
		between monotone Hurwitz numbers and Euler characteristics of Hurwitz moduli spaces \cite{EEHS}
		which seems to be very interesting, but which we do not explore in detail here.

	\section{Hurwitz-Cayley Graph}	
	\label{sec:HurwitzCayley}
	
		\subsection{Elementary features}						
		Let $\group{S}^d=\mathrm{Aut}\{1,\dots,d\}$ be the symmetric group of degree $d \in \N$, where
		we define the group law as $\pi_1\pi_2=\pi_2 \circ \pi_1$ so that 
		permutations are multiplied left to right. Let $K=\{(i\ j) \colon 1 \leq i <j \leq d\}$
		be the conjugacy class of transpositions in $\group{S}^d$. The Hurwitz-Cayley graph
		has vertex set $\group{S}^d$, with $\rho,\sigma$ adjacent
		if and only if $\rho\tau= \sigma$ for $\tau \in K$. We henceforth identify 
		$\group{S}^d$ with the Hurwitz-Cayley graph.
		
		Geodesic distance in Cayley graphs is implemented 
		by word norm. For the Hurwitz-Cayley grpah, $\mathrm{d}(\rho,\sigma)=|\rho^{-1}\sigma|$ where $|\pi|$ 
		is the minimal length of a factorization of $\pi$ into transpositions.
		The Join-Cut Lemma \cite{GJ} gives $|\pi| = d-\ell(\pi)$, where $\ell(\pi)$ is the number of 
		factors in the unique decomposition of $\pi$ as a product of disjoint cycles.
		Accordingly, the Hurwitz-Cayley graph decomposes as
						
			\begin{equation}
				\group{S}^d = \bigsqcup_{r=0}^{d-1} L_r,
			\end{equation}
			
		\noindent
		where $L_r$ is the set of permutations consisting of $d-r$ disjoint cycles. Every
		edge of the Hurwitz-Cayley graph spans consecutive levels $L_r,L_{r+1}$. 
		
		The Cayley graph of any group with respect to any symmetric generating set 
		decomposes into independent sets made up of identity-centered spheres, giving
		a graded graph structure.
		The Hurwitz-Cayley graph has the non-generic feature that
		its levels are unions of conjugacy classes: we have
		
			\begin{equation}
				L_r = \bigsqcup_{\substack{\alpha \in \Y^d \\
				\ell(\alpha)=d-r}} K_\alpha,
			\end{equation}
			
		\noindent
		where $\Y^d$ is the set of Young diagrams with exactly $d$ cells, 
		$\ell(\alpha)$ denotes the number of rows of $\alpha \in \Y^d$, and $K_\alpha$ is the 
		conjugacy class of permutations in $\group{S}^d$ of cycle type $\alpha$.
				
		\subsection{Jucys-Murphy correspondence: combinatorial form}	
		There are in general many geodesics between permutations $\rho,\sigma \in \group{S}^d$. The 
		number of geodesics $\rho \to \sigma$ equals the number of minimal transposition
		factorizations of $\rho^{-1}\sigma$, and this number was computed independently by Hurwitz
		and Cayley; see \cite{Lando:ICM}. It is 
				
			\begin{equation}
			\label{eqn:HurwitzCayley}
				{d-r \choose \alpha_1-1,\dots,\alpha_r-1} \prod_{i=1}^r \mathrm{Cay}_{\alpha_i-1},
			\end{equation}
			
		\noindent
		where $\alpha_1,\dots,\alpha_r$ are the row lengths of the Young diagram $\alpha \in \Y^d$
		encoding the cycle type of $\rho^{-1}\sigma$ and $\mathrm{Cay}_n = (n+1)^{n-1}$ is
		the Cayley number. The reasoning behind this formula is elementary. For each $i=1,\dots,r$, 
		factor cylce $i$ of $\rho^{-1}\sigma$ into $\alpha_i-1$ transpositions, the minimal 
		number required; this can be done in $\Cay_{\alpha_i-1}$ ways. The result is a string
		of $d-r$ transpositions divided into $r$ blocks, the factored cycles of $\rho^{-1}\sigma$;
		by minimality the factors in distinct blocks are disjoint. The multinomial coefficient counts shuffles of these 
		factors which maintain the relative order within each block, preserving the total product.
				
\begin{figure}
         		   \centering
        			    \includegraphics[height=8.5cm]{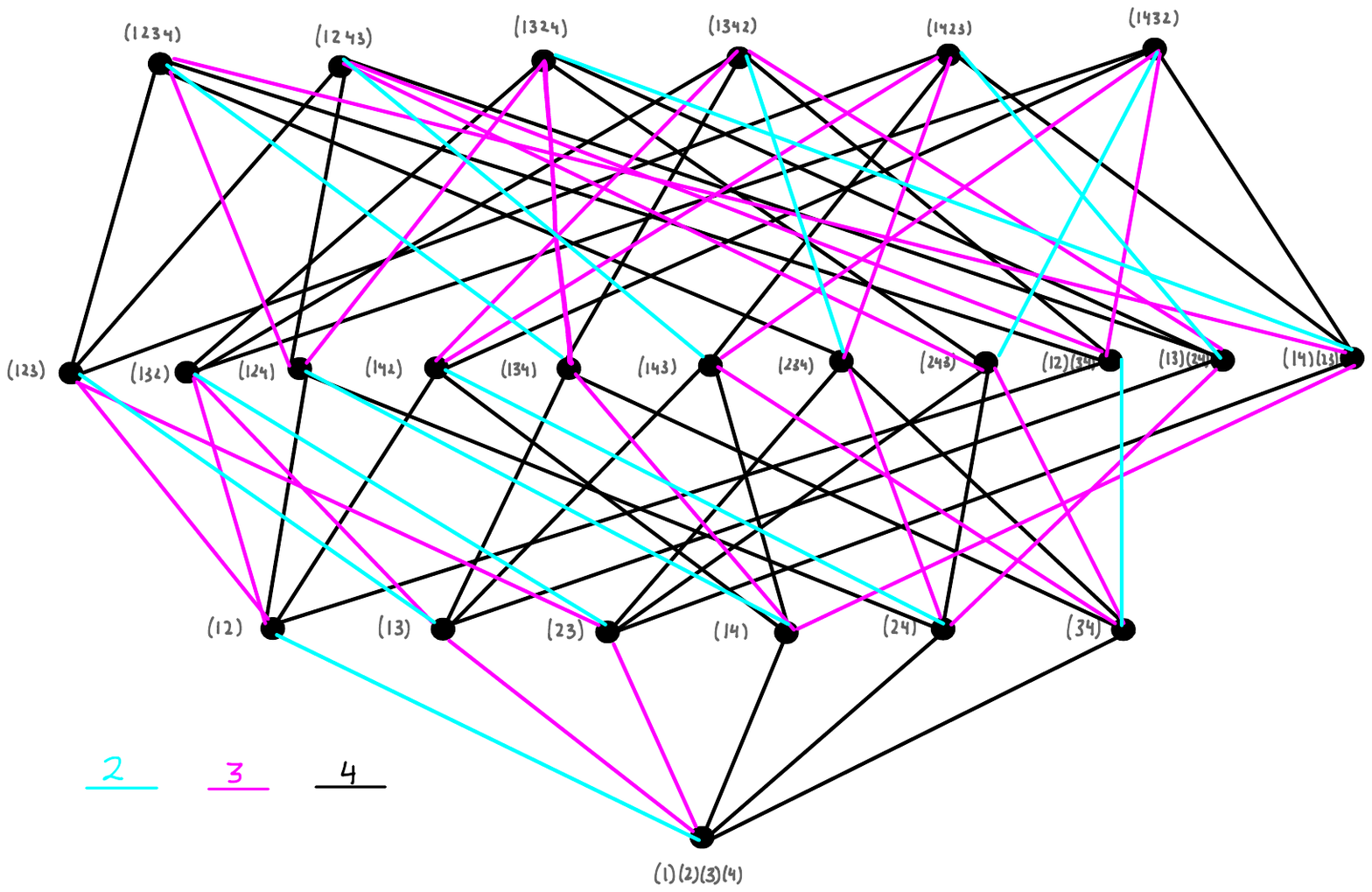}
          		 	 \caption{Jucys-Murphy labeling of $\group{S}^4$.}
           			 \label{fig:JMLabeling}
       			 \end{figure}

		A fundamental feature of the Hurwitz-Cayley graph is that it can be Euclideanized 
		by choosing a distinguished geodesic between every pair of points
		in a systematic way. This construction involves an edge labeling of the symmetric
		group introduced, implicitly and independently, by Jucys \cite{Jucys} and 
		Murphy \cite{Murphy}. Let us mark each edge corresponding to the transposition 
		$(i\ j)$ with $j$, the larger of the two numbers interchanged. Thus,
		emanating from every vertex of $\group{S}^d$ we have one $2$-edge,
		two $3$-edges, three $4$-edges, etc. See Figure \ref{fig:JMLabeling} for the case $d=4$.
		
		A walk on $\group{S}^d$ is said
		to be \emph{strictly monotone} if the labels of the edges it traverses
		form a strictly increasing sequence. Clearly the length of any such 
		walk is at most $d-1$, the number of edge labels and 
		the diameter of $\group{S}^d$. 
				
			\begin{thm}
			\label{thm:CombinatorialJM}
				There exists a unique strictly monotone walk between every 
				pair of permutations, and this walk is a geodesic.
			\end{thm}
		
		\noindent
		As with the Hurwitz-Cayley geodesic count, the reasoning behind the Jucys-Murphy 
		geodesic count is elementary. A strictly 
		monotone walk $\rho \to \sigma$ is a factorization of
		$\rho^{-1}\sigma$ into transpositions of the form
		
			\begin{equation}
				\rho^{-1}\sigma = (i_1\ j_1) \dots (i_r\ j_r), \quad j_1< \dots < j_r.
			\end{equation}
			
		\noindent
		In order to establish Theorem \ref{thm:CombinatorialJM},
		it is sufficient to establish that a cycle admits a unique strictly monotone factorization, and 
		show that the length of this factorization is minimal
		--- monotonicity kills the shuffle factor present in the unrestricted geodesic count.
		For an inductive proof of a statement equivalent to Theorem \ref{thm:CombinatorialJM},
		see \cite{DG}.
		
\begin{figure}
         		   \centering
        			    \includegraphics[height=8.5cm]{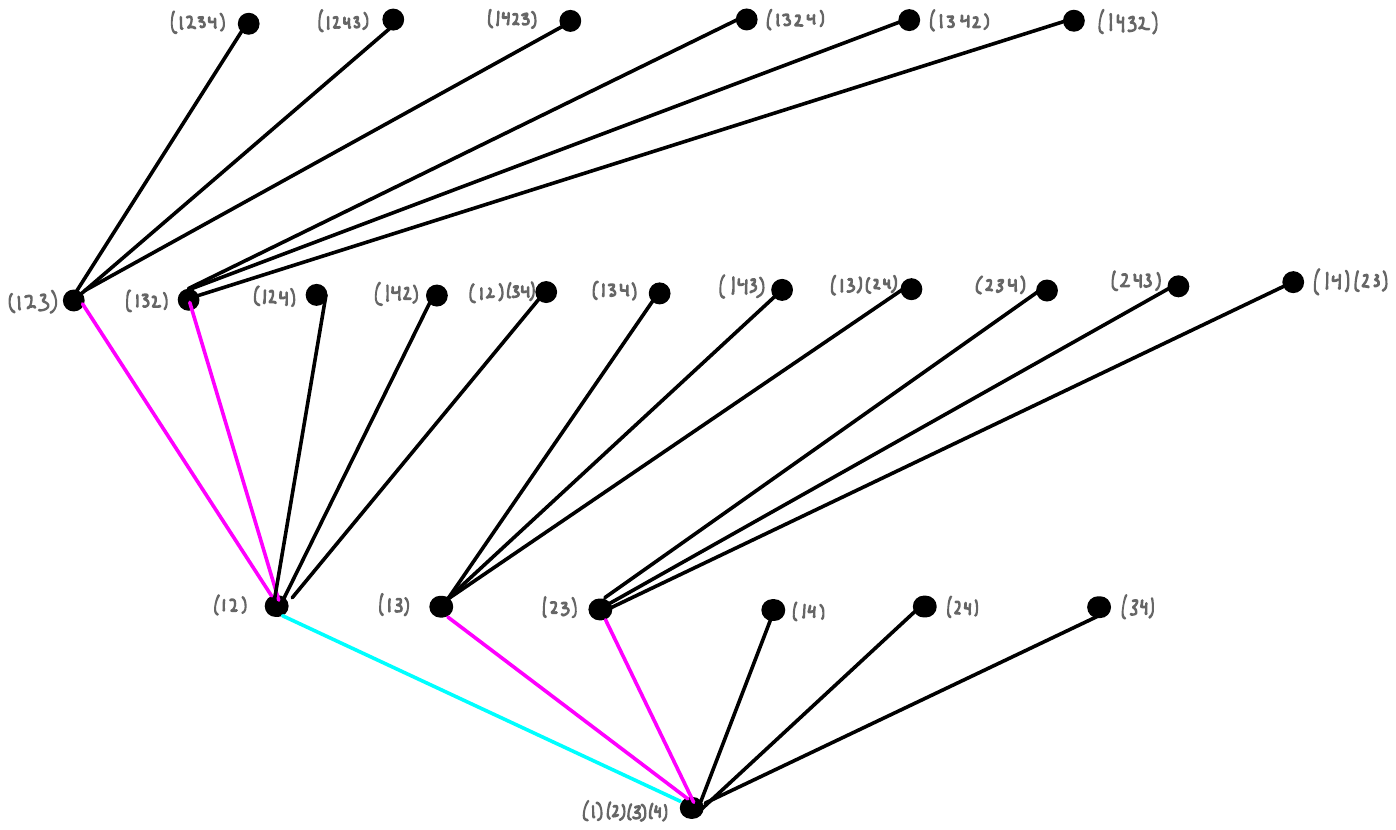}
          		 	 \caption{Jucys-Murphy correspondence in $\group{S}^4$.}
           			 \label{fig:JMCorrespondence}
       			 \end{figure}

		Theorem \ref{thm:CombinatorialJM} says that for any permutation $\rho$,
		the points of the sphere of radius $r$ centered at $\rho$ correspond
		bijectively to endpoints of strictly monotone $r$-step walks emanating
		from $\rho$. We refer to this bijection the \emph{Jucys-Murphy correspondence};
		see Figure \ref{fig:JMCorrespondence}.
		
		Instead of defining a strictly monotone walk by the condition that the labels
		of its edges are strictly increasing, we could have 
		stipulated that they strictly decrease along the walk. Since
		path reversal defines a bijection between strictly decreasing walks $\rho \to \sigma$
		and strictly increasing walks $\sigma \to \rho$, Theorem 
		\ref{thm:CombinatorialJM} implies that there is one and only one 
		strictly decreasing walk between any two vertices of the Hurwitz-Cayley graph,
		and that it is a geodesic. Theorem \ref{thm:CombinatorialJM} is thus valid 
		in both senses of the word ``monotone.''
		
		Theorem \ref{thm:CombinatorialJM} allows us to view an
		arbitrary product $\pi_1\pi_2\dots\pi_k$ of permutations
		graphically, as a concatenation of strictly monotone walks
		
			\begin{equation}
			\label{eqn:StrictConcatenation}
				\iota \to \pi_1 \to \pi_1\pi_2 \to \dots \to \pi_1\pi_2\dots\pi_k.
			\end{equation}
		
		\subsection{Jucys-Murphy correspondence: matricial form}
		Since $\group{S}^d$ is a group, $\C\group{S}^d$ is an algebra,
		the product of kets being $|\pi_1\rangle |\pi_2\rangle = |\pi_1\pi_2\rangle$. 
		The group algebra $\C\group{S}^d$ acts on itself by right multiplication: 
		every $|A \rangle \in \C\group{S}^d$ gives
		a linear operator $A \in \End\C\group{S}^d$ defined by 
		
			\begin{equation}
				A|\rho\rangle = \sum_{\pi \in \group{S}^d} \langle \pi | A \rangle |\rho\pi \rangle.
			\end{equation}
			
		\noindent
		The matrix elements of $A \in \End \C\group{S}^d$ in the 
		permutation basis are
		
			\begin{equation}
				\langle \sigma | A | \rho \rangle =  \sum_{\pi \in \group{S}^d} \langle \pi | A \rangle \langle \sigma |\rho\pi \rangle
				=\langle \rho^{-1}\sigma|A\rangle.
			\end{equation}
		
		\noindent
		In particular,
		every diagonal matrix element is equal to a common value $\langle A \rangle = \langle \iota | A\rangle$,
		the coefficient of $|\iota\rangle$ in $|A\rangle$. The linear functional $|A\rangle \mapsto \langle A \rangle$ is 
		the normalized character of the regular representation,
		
			\begin{equation}
				\langle A \rangle = \frac{1}{d!} \Tr A.
			\end{equation}
				
		The map $|A\rangle \mapsto A$ is a faithful linear representation of $\C\group{S}^d$, 
		the right regular representation. By abuse of 
		notation, if $|A\rangle \in \C\group{S}^d$ has $\{0,1\}$-coefficients 
		in the permutation basis we write $A$ for both the subset of 
		$\group{S}^d$ determined by the nonzero coefficients of $|A\rangle$,
		and for the operator on $\C\group{S}^d$ which is the image of $|A\rangle$ in the regular representation.
		Thus $K$ denotes both the conjugacy class of transpositions
		in $\group{S}^d$ and the adjacency operator on the 
		Hurwitz-Cayley graph, and $L_r$ denotes both the sphere
		of radius $r$ centered at $\iota$ in $\group{S}^d$ and the operator
		which maps $|\rho\rangle$ to the sum of all permutations $|\sigma\rangle$ on the sphere of
		radius $r$ with center $\rho$. In matrix form, the
		Jucys-Murphy correspondence is as follows.
				
			\begin{thm}
			\label{thm:MatricialJM}
				For any $\rho,\sigma \in \group{S}^d$, we have
				
					\begin{equation*}
						\langle \sigma | L_r | \rho \rangle = W^r_{<}(\rho,\sigma) 
					\end{equation*}
					
				\noindent
				where $W^r_{<}(\rho,\sigma)$ is the number of strictly 
				monotone $r$-step walks $\rho \to \sigma$ in the 
				Hurwitz-Cayley graph.
			\end{thm}
		
		Let $|\Omega_q\rangle \in \C\group{S}^d$ be the exponential distance element,
		
			\begin{equation}
				|\Omega_q \rangle = \sum_{\pi \in \group{S}^d} q^{|\pi|} |\pi\rangle = \sum_{r=0}^{d-1} q^r |L_r\rangle,
			\end{equation}
			
		\noindent
		so that $\Omega_q \in \End \C\group{S}^d$ is the exponential distance operator 
		of the Hurwitz-Cayley graph. Theorem \ref{thm:MatricialJM} immediately implies 
		the following alternative interpretation of the matrix elements of $\Omega_q$
		in the vertex basis.
		
			\begin{cor}
			For any $\rho,\sigma \in \group{S}^d$, we have
				
				\begin{equation*}
					\langle \sigma | \Omega_q | \rho \rangle = \sum_{q=0}^\infty q^rW^r_{<}(\rho,\sigma).
				\end{equation*}
			\end{cor}
					
		\subsection{Jucys-Murphy correspondence: polynomial form}
		Theorem \ref{thm:MatricialJM} says
		that 
		
			\begin{equation}
				|L_r \rangle = \sum_{2 \leq j_1 < \dots < j_r \leq d}
				\left(\sum_{i_1=1}^{j_1} |i_1\ j_1\rangle \right) \dots \left(\sum_{i_r=1}^{j_r} |i_r\ j_r\rangle \right).
			\end{equation}
			
		\noindent
		This is the statement that
		
			\begin{equation}
				|L_r \rangle = e_r(|J_1\rangle,\dots, |J_r\rangle),
			\end{equation}
			
		\noindent
		the elementary symmetric polynomial $e_r$ of degree $r$ evaluated on the 
		\emph{Jucys-Murphy elements} 
		
			\begin{equation}
			\begin{split}
				|J_1\rangle &= |0\rangle \\
				|J_2\rangle &= |1\ 2 \rangle \\
				|J_3 \rangle &= |1\ 3\rangle + |2\ 3\rangle \\\
				&\vdots \\
				|J_d \rangle &= |1\ d\rangle + \dots + |d-1\ d\rangle
			\end{split}
			\end{equation}
			
		\noindent
		in the group algebra $\C\group{S}^d$. These elements, 
		which commute with one another, are of fundamental importance
		in the representation theory of the symmetric groups. A very readable
		account of their properties has been given by Diaconis and Greene \cite{DG}.
		Their images $J_1,J_2,\dots,J_r \in \End\C\group{S}^d$ in the
		regular representation are called the \emph{Jucys-Murphy operators}.
		As a polynomial identity in $\C\group{S}^d$ or $\End \C\group{S}^d$,
		the Jucys-Murphy correspondence takes the following form.

			\begin{thm}
			\label{thm:PolynomialJM}
				For any $0 \leq r \leq d-1$, we have 
				
					\begin{equation*}
						|L_r \rangle = e_r(|J_1\rangle,\dots,|J_d\rangle)
					\end{equation*}
					
				\noindent
				in $\C\group{S}^d$, or equivalently 
				
					\begin{equation*}
						L_r = e_r(J_1,\dots,J_d)
					\end{equation*}
					
				\noindent
				in $\End \C\group{S}^d$.
			\end{thm}		
		
		For $r=1$, Theorem \ref{thm:PolynomialJM} is the obvious fact that
		the adjacency operator of the Hurwitz-Cayley graph is the sum of the 
		Jucys-Murphy operators,
		
			\begin{equation}
				K=J_1+\dots+J_d,
			\end{equation}
			
		\noindent
		which implies that the exponential adjacency operator factors as 
		
			\begin{equation}
				\Psi_t=e^{-tK} = e^{-tJ_1} \dots e^{-tJ_d}.
			\end{equation}
			
		An analogous but subtler factorization of the exponential distance operator of 
		the Hurwitz-Cayley graph is obtained by combining Theorem \ref{thm:PolynomialJM}
		with the generating series
				
			\begin{equation}
				\sum_{r=0}^\infty q^r e_r(x_1,x_2,x_3,\dots) = \prod_{i=1}^\infty (1+qx_i)
			\end{equation}
			
		\noindent
		for the elementary symmetric functions in an alphabet $x_1,x_2,x_2,\dots$ of
		commuting indeterminates \cite{GJ:Book}. 
						
			\begin{cor}
				We have 
				
					\begin{equation*}
						|\Omega_q \rangle = (|\iota\rangle+q|J_1\rangle) \dots (|\iota\rangle+q|J_d\rangle)
					\end{equation*}
				
				\noindent
				in $\C\group{S}^d$, or equivalently 
				
					\begin{equation*}
						\Omega_q = (I+qJ_1) \dots (I+qJ_d)
					\end{equation*}
					
				\noindent
				in $\End \C\group{S}^d$.
			\end{cor}
			
		We conclude that computing the eigenvalues of the exponential adjacency 
		and distance operators $\Psi_t$ and $\Omega_q$ on the Hurwitz-Cayley
		graph reduces to computing the eigenvalues of the Jucys-Murphy 
		operators $J_1,\dots,J_d \in \End\C\group{S}^d$. Before addressing this spectral problem, we consider a 
		further implication of the Jucys-Murphy correspondence.
		
		\subsection{Weakly monotone walks}
		Combining the polynomial form of the Jucys-Murphy correspondence
		with Newton's theorem on symmetric polynomials,
		every symmetric polynomial function of the Jucys-Murphy elements
		$|J_1\rangle,\dots,|J_d\rangle \in \C\group{S}^d$
		is a polynomial function of the levels $|L_r\rangle = e_r(|J_1\rangle,\dots,|J_d\rangle)$
		of the Hurwitz-Cayley graph. As explained above, each level is a sum of 
		conjugacy classes,
			
			\begin{equation}
			\label{eqn:LevelClassSum}
				|L_r\rangle=\sum_{\substack{\alpha \in \Y^d \\ \ell(\alpha)=d-r}} |K_\alpha\rangle, \quad 0 \leq r \leq d-1.
			\end{equation}
			
		\noindent
		Consequently, every symmetric polynomial function $f(|J_1\rangle,\dots,|J_d\rangle)$
		is a linear combination of conjugacy classes: we have
				
			\begin{equation}
				f(|J_1\rangle,\dots,|J_d\rangle) = \sum_{\alpha \in \Y^d} c_\alpha(f) |K_\alpha \rangle
			\end{equation}
			
		\noindent
		for some coefficients $c_\alpha(f)$, which are integers if $f$ is an integral
		polynomial in the $e_r$'s. 
				
		A particularly interesting case is that of complete symmetric polynomials in Jucys-Murphy
		elements,
		
			\begin{equation}
				|M_r\rangle = h_r(|J_1\rangle,\dots,|J_d\rangle) = \sum_{2 \leq j_1 \leq \dots \leq j_r \leq d}
				|J_{j_1}\rangle \dots |J_{j_r}\rangle.
			\end{equation}
			
		\noindent
		We say that a walk on the Hurwitz-Cayley graph is \emph{weakly monotone} if the labels
		of the edges it traverses form a weakly increasing sequence, and write 
		$W^r_{\leq}(\rho,\sigma)$ for the number of weakly monotone $r$-step walks 
		$\rho \to \sigma$ between given permutations $\rho,\sigma \in \group{S}^d$. 
		By definition of the Jucys-Murphy elements and the complete symmetric
		polynomials, we have that
		
			\begin{equation}
				\langle \pi | M_r \rangle = W^r_{\leq}(\iota,\pi)
			\end{equation}
			
		\noindent
		is the number of weakly monotone $r$-step walks $\iota \to \pi$. 		
		Equivalently, for $M_r \in \End \C\group{S}^d$ the image of $|M_r\rangle$
		in the regular representation, we have the matrix formula
		
			\begin{equation}
				\langle \sigma | M_r | \rho \rangle = W^r_{\leq}(\rho,\sigma),
			\end{equation}
			
		\noindent
		the number of $r$-step weakly monotone walks $\rho \to \sigma$ in the
		Hurwitz-Cayley graph. 
		
		The fact that $|M_r\rangle$ is a central element in $\C\group{S}^d$
		implies that $W^r_{\leq}(\iota,\pi)$ depends only on the conjugacy class
		of $\pi \in \group{S}^d$, which is equivalent to the statement that
		$W^r_{\leq}(\rho,\sigma)$ depends only on the cycle type 
		of $\rho^{-1}\sigma$. In particular, we have
				
			\begin{equation}
				W^r_{\leq}(\sigma,\rho) = W^r_{\leq}(\iota,\sigma^{-1}\rho) =
				W^r_{\leq}(\iota,\rho^{-1}\sigma) = W^r_{\leq}(\rho,\sigma),
			\end{equation}
			
		\noindent
		which by path reversal implies that 
		
			\begin{equation}
				W^r_{\leq}(\rho,\sigma) = W^r_{\geq}(\rho,\sigma),
			\end{equation}
			
		\noindent
		where $W^r_{\geq}(\rho,\sigma)$ is the number of $r$-step 
		walks $\rho \to \sigma$ on the Hurwitz-Cayley graph whose
		step labels form a weakly decreasing sequence. Thus
		at the enumerative level there is no difference
		between decreasing and increasing walks in the 
		Hurwitz-Cayley graph, whether strictly or weakly.
		
		We have seen above that the total number of geodesics $\rho \to \sigma$
		in the Hurwitz-Cayley graph is given by a product of Cayley numbers along
		the cycles of $\rho^{-1}\sigma$ times a shuffle factor. We have also seen
		that the number of strictly monotone geodesics $\rho \to \sigma$ is exactly 
		one. The number of weakly monotone geodesics $\rho \to \sigma$, i.e. the
		number of minimal weakly monotone factorizations of $\rho^{-1}\sigma$ into
		transpositions, must lie between these two extremes.
		
			\begin{thm}[\cite{Novak:BCP}]
			\label{thm:MonotoneGeo}
				The number of weakly monotone geodesics $\rho \to \sigma$
				in the Hurwitz-Cayley graph $\group{S}^d$ is
				
					\begin{equation*}
						\prod_{i=1}^{\ell(\alpha)} \Cat_{\alpha_i-1},
					\end{equation*}
					
				\noindent
				where $\alpha \in \Y^d$ is the cycle type of $\rho^{-1}\sigma$
				and $\Cat_n=\frac{1}{n+1}{2n \choose n}$ is the Catalan number.
			\end{thm}
		
		The proof of Theorem \ref{thm:MonotoneGeo} is again elementary, reducing
		to counting minimal monotone factorizations of a cycle; that this is a Catalan
		number was first discovered by Gewurz and Merola \cite{GM}, and the lack
		of shuffle factor is again due to monotonicity. A more general result counting
		weakly monotone factorizations with a given ``signature,'' corresponding to evaluation
		of monomial symmetric polynomials in Jucys-Murphy elements, was proved
		in \cite{MatNov}. This result can in turn be used to develop an enumerative theory of weighted
		walks in the Hurwitz-Cayley graph \cite{GH}, which turns out to be very rich \cite{ACEH}.
									
		\subsection{Fourier transform}
		As with the group algebra of any finite group \cite{Serre}, we have an orthogonal decomposition 
				
			\begin{equation}
			\label{eqn:FourierIsomorphism}
				\C\group{S}^d = \bigoplus_{\lambda \in \Y^d} (\dim \V^\lambda) \V^\lambda 
			\end{equation}
			
		\noindent
		where $\V^\lambda$ is an enumeration of the irreducible representations of $\C\group{S}^d$.
		By Schur's Lemma, if $| A \rangle \in \C\group{S}^d$
		is a linear combination of conjugacy classes then $A \in \End \C\group{S}^d$ acts
		in $\V^\lambda$ as a scalar operator whose eigenvalue we denote $\hat{A}(\lambda)$. 
		This associates to every central element $|A\rangle \in \C\group{S}^d$ a function
		$\hat{A} \colon \Y^d \to \C$, the \emph{Fourier transform} of $|A \rangle$. 
		The map $|A\rangle \to \hat{A}$ is an isometric algebra isomorphism from the center of $\C\group{S}^d$ to 
		the function algebra $L^2(\Y^d)$ --- we have
		
			\begin{equation}
				\langle K_\alpha | K_\beta \rangle = \sum_{\lambda \in \Y^d} \frac{(\dim \V^\lambda)^2}{d!}
				\hat{K}_\alpha(\lambda) \hat{K}_\beta(\lambda) =\delta_{\alpha\beta} |K_\alpha|.
			\end{equation}
			
		\noindent
		This gives the Plancherel expectation formula 
						
			\begin{equation}
				\langle A \rangle = \frac{1}{d!} \Tr A = \sum_{\lambda \in \Y^d} 
				\frac{(\dim \V^\lambda)^2}{d!} \hat{A}(\lambda).
			\end{equation}
			
		\noindent
		for central elements $|A\rangle \in \C\group{S}^d$.
					
		A remarkable formula for computing the Fourier transform of a symmetric polynomial
		function in Jucys-Murphy elements was discovered by Jucys and Murphy; see \cite{DG}.
		Recall that the
		content $c(\Box)$ of a cell $\Box$ in a Young diagram $\lambda$ is defined to 
		be its column index minus its row index. Thus, filling the cells of $\lambda$ with
		their contents produces a jagged $\lambda$-shaped corner of the infinite 
		Toeplitz matrix $[j-i]_{i,j=1}^\infty$. 
		
			\begin{thm}
			\label{thm:SpectralJM}
				If $|A\rangle \in \C\group{S}^d$ is such that $|A\rangle = f(|J_1\rangle,\dots,|J_d\rangle)$
				for a symmetric polynomial $f$, then the Fourier transform of $|A\rangle$ is
				
					\begin{equation*}
						\hat{A}(\lambda) = f(c(\Box) \colon \Box \in \lambda),
					\end{equation*}
					
				\noindent
				the evaluation of $f$ on the multiset of conents of $\lambda$. 
				Equivalently, for each $\lambda \in \Y^d$ the
				eigenvalue of $A = f(J_1,\dots,J_d) \in \End \C\group{S}^d$
				acting in $\V^\lambda$ is $f(c(\Box) \colon \Box \in \lambda)$.
			\end{thm}
		
		According to Theorem \ref{thm:SpectralJM}, the transposition class
		$|K\rangle = |J_1\rangle + \dots + |J_d\rangle$ has Fourier transform
		
			\begin{equation}
				\hat{K}(\lambda) = \sum_{\Box \in \lambda} c(\Box).
			\end{equation}
			
		\noindent
		The spectrum of the adjacency operator of the Hurwitz-Cayley graph consists of these
		numbers as $\lambda$ ranges over $\Y^d$,
		with the multiplicity of $\hat{K}(\lambda)$ being $\dim \V^\lambda$.
		The eigenvalues of the exponential
		adjacency operator are thus
		
			\begin{equation}
			\label{eqn:ExponentialAdjacencyEigenvalues}
				\hat{\Psi}_t(\lambda) = \prod_{\Box \in \lambda} e^{-tc(\Box)}, \quad \lambda \in \Y^d.
			\end{equation}
		
		To find the eigenvalues of the exponential distance operator $\Omega_q$ of the 
		Hurwitz-Cayley graph, we need to compute the Fourier of every level of the graph.
		Theorem \ref{thm:SpectralJM} gives the Fourier transform of
		$|L_r\rangle = e_r(|J_1\rangle,\dots,|J_d\rangle)$ as
		
			\begin{equation}
				\hat{L}_r(\lambda) = e_r(c(\Box) \colon \Box \in \lambda).
			\end{equation}
			
		\noindent
		The eigenvalues of the exponential adjacency operator 
		
			\begin{equation}
				\Omega_q = \sum_{r=0}^{d-1} q^r L_r
			\end{equation}
			
		\noindent
		are thus
		
			\begin{equation}
			\label{eqn:ExponentialDistanceEigenvalues}
				\hat{\Omega}_q(\lambda) = \prod_{\Box \in \lambda} (1+qc(\Box)), \quad
				\lambda \in \Y^d, 
			\end{equation}
	
		\noindent
		and
		
			\begin{equation}
				\det \Omega_q = \prod_{\lambda \in \Y^d} \hat{\Omega}_q(\lambda)^{\dim \V^\lambda} = 
				\prod_{\lambda \in \Y^d} \prod_{\Box \in \lambda} (1+qc(\Box))^{\dim \V^\lambda}.
			\end{equation}
			
		\noindent
		One can give an analogous formula for any coefficient
		of the characteristic polynomial of $\Omega_q$, but we will not do so here.
		The determinant formula is enough to see the following.
			
			\begin{thm}
				The exponential distance operator $\Omega_q$ of the 
				Hurwitz-Cayley graph $\group{S}^d$ is singular if and only if $d>1$ and
			 	$q$ or $-q$ is one of the unit fractions
					
					\begin{equation*}
						\frac{1}{d-1}, \frac{1}{d-2},\dots,\frac{1}{1}.
					\end{equation*}
			\end{thm}
			
		\subsection{The operator $\Omega_q^{-1}$}
		For $q$ away from the above unit fractions, the exponential distance operator 
		$\Omega_q$ of the Hurwitz-Cayley graph is invertible, and 
		its inverse acts in $\V^\lambda$ as multiplication by 
		
			\begin{equation}
				\hat{\Omega}_q^{-1}(\lambda) = \prod_{\Box\in \lambda} \frac{1}{1+qc(\Box)}.
			\end{equation}
			
		\noindent
		For $|q| < \frac{1}{d-1}$, this is the absolutely convergent series 
		
			\begin{equation}
				\hat{\Omega}_q^{-1}(\lambda) = \sum_{r=0}^\infty (-q)^r h_r(c(\Box) \colon \Box \in \lambda)
			 = \sum_{r=0}^\infty (-q)^r \hat{M}_r(\lambda),
			\end{equation}
			
		\noindent
		where $|M_r\rangle=h_r(|J_1\rangle,\dots,|J_d\rangle)$. 
		Thus, for $|q| < \frac{1}{d-1}$, we have
		
			\begin{equation}
				\Omega_q^{-1} = \sum_{r=0}^\infty (-q)^r M_r
			\end{equation}
			
		\noindent
		in $\End \C\group{S}^d$, where $M_r=h_r(J_1,\dots,J_d)$ is 
		the complete symmetric polynomial of degree $r$ in the Jucys-Murphy 
		operators. The matrix elements of the inverse exponential distance 
		operator are therefore
		
			\begin{equation}
				\langle \sigma | \Omega_q^{-1} | \rho\rangle = \sum_{r=0}^\infty (-q)^r 
				\langle \sigma | M_r | \rho \rangle,
			\end{equation}
			
		\noindent
		and we get a combinatorial interpretation for the matrix elements of $\Omega_q^{-1}$ 
		as generating functions for weakly monotone walks in the Hurwitz-Cayley graph.
		
			\begin{thm}
				For $|q|<\frac{1}{d-1}$, the matrix elements of the inverse exponential distance operator
				on $\Omega_q$ are absolutely convergent generating functions 
												
					\begin{equation*}
						\langle \sigma | \Omega_q^{-1} | \rho \rangle = \sum_{r=0}^\infty (-q)^rW^r_{\leq}(\rho,\sigma)
					\end{equation*}
					
				\noindent
				for weakly monotone walks in the Hurwitz-Cayley graph $\group{S}^d$, 
			\end{thm}

		\subsection{Summary}
		The exponential adjacency and distance matrices $\Psi_t$ and $\Omega_q$
		of the Hurwitz-Cayley graph are generating functions for walks in 
		$\group{S}^d$. The matrix elements of
		the exponential adjacency matrix are 
		
			\begin{equation}
				\langle \sigma | \Psi_t | \rho \rangle = \sum_{r=0}^\infty \frac{(-t)^r}{r!} W^r(\rho,\sigma),
			\end{equation}
			
		\noindent
		where $W^r(\rho,\sigma)$ is the total number of walks $\rho \to \sigma$. This is a general
		fact about graphs. The Hurwitz-Cayley graph has the special feature that
				
			\begin{equation}
				\langle \sigma | \Omega_q | \rho\rangle = \sum_{r=0}^\infty q^r W^r_{<}(\rho,\sigma)
				\quad\text{and}\quad
				\langle \sigma | \Omega_q^{-1} | \rho\rangle = \sum_{r=0}^\infty (-q)^r W^r_{\leq}(\rho,\sigma)
			\end{equation}
			
		\noindent
		are generating functions enumerating strictly and weakly monotone
		walks $\rho \to \sigma$ in $\group{S}^d$.
		The first series consists of a single nonzero term, namely $q^{\mathrm{d}(\rho,\sigma)}$. 
		The second series has infinitely many nonzero terms and requires
		$|q|<\frac{1}{d-1}$ to ensure absolute convergence.
				
\section{Yang-Mills Theory}
\label{sec:YangMills}

\subsection{Partition functions}
The partition function of Yang-Mills theory on a compact two-dimensional orientable spacetime of area $t$ with $m$ holes and $n$ handles
admits a dual representation in terms of the gauge group, $\group{U}_N$. 
The elegant formula reads \cite{CMR1}

	\begin{equation}
	\label{eqn:PartitionFunction}
		\mathcal{Z}_N = \sum_{\lambda_1 \geq \dots \geq \lambda_N} \frac{S_\lambda(U_1,\dots,U_m)}{(\dim \W_N^\lambda)^{m+2n-2}}
		e^{-\frac{t}{2N} \hat{C}(\lambda)}.
	\end{equation}
	
\noindent
The sum is over integer vectors $\lambda=(\lambda_1,\dots,\lambda_N)$ 
with weakly decreasing coordinates, each labeling an irreducible representation 
$\W_N^\lambda$ of $\group{U}_N$. The function $S_\lambda \colon \group{U}_N^m \to \C$ 
is the character of $(\W_N^\lambda)^{\otimes m}$. The unitary matrices
$U_1,\dots,U_m$ represent boundary holonomies: setting any one of them 
to the identity cancels a factor of $\dim \W_N^\lambda$ and plugs a hole. The exponential
factor is a discrete Gaussian weight: $\hat{C}(\lambda)$ is the $\W_N^\lambda$-eigenvalue of the quadratic Casimir $C$,
a central element in the universal enveloping algebra of $\mathfrak{u}_N$ corresponding to the Laplacian on $\group{U}_N$. 
If spacetime is a cylinder, then $\mathcal{Z}_N$ is the $\group{U}_N$
heat kernel \cite{KW,Zelditch}.

It was perceived by Gross that \eqref{eqn:PartitionFunction} provides a promising 
path to a string description of $YM_2$ at large $N$, which ought to be Hurwitz theory \cite{Gross}. 
This insight was fully developed by Gross and Taylor \cite{GT1,GT2}, who
argued that as $N \to \infty$ the partition function $\mathcal{Z}_N$ factorizes
into two copies of a \emph{chiral partition function}, $Z_N$, obtained by restricting the sum 
\eqref{eqn:PartitionFunction} to nonnegative $\lambda$.
The chiral partition function may be viewed as a sum over the set $\Y_N$ of Young diagrams 
with at most $N$ rows, and thus decomposes as

	\begin{equation}
	\label{eqn:ChiralStratified}
		Z_N = 1 + \sum_{d=1}^\infty Z_N^d,
	\end{equation}
	
\noindent
where the \emph{microchiral partition function}

	\begin{equation}
	\label{eqn:Microchiral}
		Z_N^d = \sum_{\lambda \in \Y_N^d}
		 \frac{S_\lambda(U_1,\dots,U_m)}{(\dim \W_N^\lambda)^{m+2n-2}}
		e^{-\frac{t}{2N} \hat{C}(\lambda)}
	\end{equation}
	
\noindent
is a sum over the finite set $\Y_N^d$ of Young diagrams with at most $N$ rows
and exactly $d$ cells. The microchiral partition function \eqref{eqn:Microchiral} is the core of $YM_2$. 
What we seek is a $1/N$ expansion of $Z_N^d$ whose coefficients count maps 
of worldsheets into spacetime. Once this string signature has been found, we
sum on $d$ to recover the chiral partition function $Z_N$, and then square
to recover the full partition function $\mathcal{Z}_N$.

\subsection{Gross-Taylor formula}
Gross and Taylor's fundamental insight \cite{GT1,GT2} is that $Z_N^d$ can be presented entirely in terms
of the representation theory of the symmetric group $\group{S}^d$.
This is accomplished in three steps.

	\subsubsection{Laplacian swap}
	The first observation of Gross and Taylor \cite{GT1} is that 
	
		\begin{equation}
		\label{eqn:Casimir}
			\hat{C}(\lambda) = dN + 2\hat{K}(\lambda), \quad \lambda \in \Y_N^d,
		\end{equation}
	
	\noindent
	where $\hat{K}(\lambda)$ is the Fourier transform of the transposition 
	class $|K \rangle \in \C\group{S}^d$.
	Thus \eqref{eqn:Casimir} effectively trades the Laplacian on $\group{U}_N$ for the Laplacian on 
	$\group{S}^d$, yielding
	
		\begin{equation}
			\label{eqn:StepOneMicrochiral}
				Z_N^d = z^d \sum_{\lambda \in \Y_N^d}  \frac{S_\lambda(U_1,\dots,U_m)}{(\dim \W_N^\lambda)^{m+2n-2}} 
				\hat{\Psi}_\frac{t}{N}(\lambda),
		\end{equation}
		
	\noindent
	where $z=e^{-t/2}$ and $\hat{\Psi}_t(\lambda) = e^{-t\hat{K}(\lambda)}$.

	\subsubsection{Character swap}
	According to Frobenius \cite{Zagier}, the character $s_\lambda$ of $\W_N^\lambda$
	is given by
		
		\begin{equation}
		\label{eqn:Frobenius}
			s_\lambda(U) = \frac{\dim \V^\lambda}{d!} \sum_{\alpha \in \Y^d} p_\alpha(U) \hat{K}_\alpha(\lambda),
		\end{equation}
		
	\noindent
	where $\hat{K}_\alpha(\lambda)$ is the Fourier transform of the conjugacy class $|K_\alpha\rangle \in \C\group{S}^d$
	and 
	
		\begin{equation}
			p_\alpha(U) = \prod_{i=1}^{\ell(\alpha)} \Tr U^{\alpha_i}
		\end{equation}
		
	\noindent
	is the corresponding trace invariant.
	We can thus replace $S_\lambda = s_\lambda \otimes \dots \otimes s_\lambda$
	with
	
		\begin{equation}
			S_\lambda = \left( \frac{\dim \V^\lambda}{d!} \right)^m \sum_{\alpha^1,\dots,\alpha^m \in \Y^d} 
			P_{\alpha^1\dots\alpha^m} \hat{K}_{\alpha^1 \dots \alpha^m}(\lambda),
		\end{equation}
		
	\noindent
	where 
	
		\begin{equation}
			P_{\alpha^1\dots\alpha^m}(U_1,\dots,U_m) = p_{\alpha^1}(U_1) \dots p_{\alpha^m}(U_m)
		\end{equation}
		
	\noindent
	is a product of trace invariants on $\group{U}_N$ and 
	
		\begin{equation}
			\hat{K}_{\alpha^1\dots\alpha^m}(\lambda) = \hat{K}_{\alpha^1}(\lambda) \dots \hat{K}_{\alpha^m}(\lambda)
		\end{equation}
		
	\noindent
	is a product of central characters of $\group{S}^d$.
	This gives the microchiral partition function as
	
		\begin{equation}
			\label{eqn:StepOneMicrochiral}
				Z_N^d = z^d \sum_{\alpha^1,\dots,\alpha^m \in \Y^d} P_{\alpha^1\dots\alpha^m}
				\sum_{\lambda \in \Y_N^d} \left( \frac{\dim \V^\lambda}{d!} \right)^m
				\frac{\hat{K}_{\alpha^1 \dots \alpha^m}(\lambda)}{(\dim \W_N^\lambda)^{m+2n-2}}
				\hat{\Psi}_\frac{t}{N}(\lambda),
		\end{equation}
		
	\noindent
	where we now view $Z_N^d$ as a function on $\group{U}_N^m$.

	\subsubsection{Dimension swap}
	Finally, we eliminate $\dim \W_N^\lambda$
	using the proportionality
	
		\begin{equation}
		\label{eqn:DimensionRatio}
			\frac{\dim \V^\lambda}{\dim \W_N^\lambda} = \frac{d!}{(N)_\lambda}, \quad \lambda \in \Y_N^d,
		\end{equation}
		
	\noindent
	where

		\begin{equation}
			(x)_\lambda = \prod_{\Box \in \lambda} (x+c(\Box))
		\end{equation}
	
	\noindent
	is the generalized Pochammer symbol, which 
	after renormalization is an eigenvalue
	
		\begin{equation}
			\hat{\Omega}_q(\lambda) = \prod_{\Box \in \lambda} (1+qc(\Box))
		\end{equation}
		
	\noindent
	of the exponential distance operator $\Omega_q$
	of the Hurwitz-Cayley graph $\group{S}^d$. We thus have
	
		\begin{equation}
		\label{eqn:DimensionRatioRenormalized}
			\frac{\dim \V^\lambda}{\dim \W_N^\lambda} = \frac{d!}{N^d} \hat{\Omega}_\frac{1}{N}^{-1}(\lambda),
		\end{equation}
	
	\noindent	
	where $\hat{\Omega}^{-1}_{\frac{1}{N}}(\lambda)$ is well-defined for any 
	Young diagram with at most $N$ rows. The microchiral partition function 
	thus assumes its final form.
	
		\begin{thm}[Gross-Taylor formula]
		\label{thm:GrossTaylor}
		The microchiral partition function of $YM_2$
		on a compact orientable surface of area $t$ with $m$ holes and $n$ handles is
		
			\begin{equation*}
				Z_N^d = \left( zN^{2-2n-m} \right)^d
				\sum_{\alpha^1,\dots,\alpha^m \in \Y^d}
				P_{\alpha^1 \dots \alpha^m}
				\sum_{\lambda \in \Y_N^d} \left(\frac{\dim \V^\lambda}{d!}\right)^{2-2n}
				\hat{K}_{\alpha^1 \dots \alpha^m}(\lambda) \hat{\Psi}_\frac{t}{N}(\lambda)\hat{\Omega}_\frac{1}{N}^{2-2n-m}(\lambda).
			\end{equation*}
	\end{thm}
		
\subsection{Great expectations}
	The Gross-Taylor formula reveals that the existence of a $1/N$ expansion of
	$Z_N^d$, and the Hurwitz-theoretic interpretation of its coefficients, is a somewhat
	subtle issue. The internal sum in Theorem \ref{thm:GrossTaylor} can be written
	
		\begin{equation}
		\label{eqn:InternalSum}
			\frac{1}{d!}\sum_{\lambda \in \Y_N^d} \frac{(\dim \V^\lambda)^2}{d!}
				\hat{H}^n(\lambda)\hat{K}_{\alpha^1 \dots \alpha^m}(\lambda) 
				\hat{\Psi}_\frac{t}{N}(\lambda)\hat{\Omega}_\frac{1}{N}^{2-2n-m}(\lambda),
		\end{equation}
		
	\noindent
	where $|H \rangle \in \C\group{S}^d$ is the central element whose Fourier transform is 
	$\hat{H}(\lambda)=\left(\frac{d!}{\dim \V^\lambda}\right)^2$. This element is 
	the commutator sum \cite{Zagier}
	
		\begin{equation}
			|H\rangle = \sum_{\rho,\sigma \in \group{S}^d} \rho^{-1}\sigma^{-1}\rho\sigma.
		\end{equation}

	\noindent
	In the \emph{stable range}, where $1 \leq d \leq N$, we have $\Y_N^d = \Y^d$
	and the sum \eqref{eqn:InternalSum} runs over all irreducible representations of $\C\group{S}^d$
	and is a Plancherel expectation.
						
		\begin{definition}
		\label{def:GreatExpectation}
			For any integers $1 \leq d \leq N$ and $m,n \in \N_0$, the corresponding
			\emph{Gross-Taylor expectation} is the function on $(\Y^d)^m$ defined by
			
				\begin{equation*}
					\mathcal{E}_N^d(\alpha_1,\dots,\alpha_m;n) = 
					\langle H^n K_{\alpha^1 \dots \alpha^m}\Psi_\frac{t}{N}\Omega_\frac{1}{N}^{2-2n-m} \rangle.
				\end{equation*}
		\end{definition}
	
	The Gross-Taylor expectation $\mathcal{E}_N^d(\alpha_1,\dots,\alpha_m;n)$ admits
	an absolutely convergent series expansion in powers of $1/N$ determined by the 
	$1/N$ expansions of the exponential adjacency operator

		\begin{equation}
			\Psi_\frac{t}{N} = \sum_{r=0}^\infty \frac{(-t)^r}{r!N^r} K^r
		\end{equation}
		
	\noindent
	and exponential distance operator and its inverse,
	
		\begin{equation}
			\Omega_\frac{1}{N} = \sum_{r=0}^{d-1} \frac{1}{N^r} L_r
			\quad\text{and}\quad
			\Omega_\frac{1}{N}^{-1} = \sum_{r=0}^\infty \frac{(-1)^r}{N^r} M_r,
		\end{equation}
		
	\noindent
	of the Hurwitz-Cayley graph $\group{S}^d$. The coefficients 
	
		\begin{equation}
			L_r=e_r(J_1,\dots,J_d)
			\quad\text{and}\quad
			M_r=h_r(J_1,\dots,J_d)
		\end{equation}
		
	\noindent 
	are elementary and complete symmetric polynomial functions of the
	Jucys-Murphy operators,
	as discussed in Section \ref{sec:HurwitzCayley}. The 
	coefficients of the $1/N$ expansion of any Gross-Taylor expectation
	are expectations of words in the commuting operators $H,K_\alpha,L_r,M_s$.
	
	In the \emph{unstable range}, where $d >N$, two problems arise. The first
	is that $\Y_N^d$ is a proper subset of $\Y^d$ and the internal sum in the 
	Gross-Taylor formula is not a Plancherel expectation. The second is that,
	if the number of holes in spacetime exceeds its Euler characteristic, 
	the exponential distance operator $\Omega_\frac{1}{N}$ appears to a negative 
	power and the resulting $1/N$ expansion is divergent.
	
	The standard approach to these issues in $YM_2$ is to ignore them and say that we take
	$N \to \infty$. This is perfectly fine at the microchiral level, where $d \in \N$ is 
	fixed and we only need $N \geq d$. However, if we wish to obtain a $1/N$ expansion 
	for the microchiral partition function $Z_N = 1+\sum_d Z_N^d$ by summing the 
	$1/N$ expansion of $Z_N^d$, an interchange of limits is required. What is actually 
	done in the literature is to deceptively replace the numerical quantity with an 
	indeterminate deceptively named $1/N$. This leads to a formal power series,
	the \emph{chiral Gross-Taylor series} \cite{CMR1,CMR2}, and understanding when this operation
	produces a quantitatively correct $N \to \infty$ asymptotic expansion of $Z_N$
	is an apparent gap in the $YM_2$ literature.
	
	\subsection{Lattice Yang-Mills}
		There is an interesting analogy with Yang-Mills on the 
		lattice, where the role of the chiral partition function is roughly played by the 
		partition function of  the Bars-Green/Br\'ezin-Gross-Witten/Wadia unitary matrix model 
		\cite{BarsGreen,BrezinGross,GrossWitten,Wadia},
			
		\begin{equation}
			W_N = \int_{\group{U}_N} e^{\sqrt{z}N\Tr (AU + BU^*)} \mathrm{d}U
		\end{equation}
		
		\noindent
		which plays a basic role in lattice gauge
		theory of any dimension. By character expansion \cite{Samuel},
	
		\begin{equation}
			W_N = 1 + \sum_{d=1}^\infty W_N^d,
		\end{equation}
		
		\noindent
		where 
	
		\begin{equation}
		\label{eqn:Weingarten}
			W_N^d = \frac{(zN)^d}{d!} \sum_{\alpha \in \Y^d} p_\alpha(AB) \sum_{\lambda \in \Y_N^d} 
			\frac{(\dim \V^\lambda)^2}{d!} K_\alpha(\lambda) \Omega_\frac{1}{N}^{-1}(\lambda)
		\end{equation}
		
	\noindent
	is the lattice counterpart of the microchiral partition function $Z_N^d$ in $YM_2$. The internal sum over $\lambda$
	in \eqref{eqn:Weingarten} has come to be known as the \emph{Weingarten function}; see \cite{CMN}.
	
	In the one-plaquette model, $\Omega$ is the sole source of $1/N$ --- everything is determined by the spectral theory 
	of the exponential distance operator on the Hurwitz-Cayley graph. The divergence of the $1/N$ expansion, 
	which is a consequence of the instability of the zeros of $\det \Omega_q$, was understood early on 
	to be a fact of life on the lattice \cite{DtH}, and diagrammatic expansions are
	very complicated \cite{BB,Creutz,OZ,Weingarten} if one works directly with $W_N^d$ as a sum of link integrals,
	
		\begin{equation}
		\label{eqn:LinkIntegral}
			W_N^d = \sum A_{j_1i_1} \dots A_{j_di_d} B_{i'_1j'_1} \dots B_{i'_dj'_d}
			\int_{\group{U}_N} U_{i_1j_1} \dots U_{i_dj_d} \overline{U_{i'_1j'_1} \dots U_{i'_dj'_d}}\ \mathrm{d}U.
		\end{equation}
		
	\noindent
	The use of monotone walks in the Hurwitz-Cayley graph as ``dual'' Feynman diagrams for 
	unitary matrix integrals was developed, in the language of permutation factorizations, in \cite{MatNov,Novak:BCP}.
	It can be quite useful in physical applications, see e.g. \cite{BK}.
	
	\subsection{Summary}
	We have derived the Gross-Taylor formula, Theorem \ref{thm:GrossTaylor}, which presents
	the microchiral partition function $Z_N^d$ of $YM_2$ on any compact orientable two-dimensional space
	time in terms of the representation theory of the symmetric group $\group{S}^d$. 
	A key feature of the formula is the appearance of the eigenvalues of the exponential
	adjacency and distance operators $\Psi_t$ and $\Omega_q$ of the Hurwitz-Cayley 
	graph, with $t \rightsquigarrow t/N$ and $q \rightsquigarrow 1/N$.
	Iin the stable range, $1 \leq d \leq N$, this combines with the power
	series expansions of these operators in $t$ and $q$ to yield an absolutely 
	convergent expansion $Z_N^d$ in powers of $1/N$, whose coefficients 
	are joint Plancherel moments of four basic commuting operators $H,K_\alpha,L_r,M_s$.
	These basic expectations admit combinatorial interpretations as counting
	trajectories in the Hurwitz-Cayley graph. 
	
	In the remainder of the paper, we calculate the $1/N$ expansion 
	of the microchiral partition function of $YM_2$ for representative choices
	of compact orientable spacetimes, and compute the corresponding chiral 
	Gross-Taylor series $Z$ in each case. We consider the large $N$ asymptotic expansion of 
	$Z_N$ predicted by the Gross-Taylor series, but do not address its quantitative 
	correctness, though this is a very interesting topic \cite{BouSch,Zelditch}.			
		
\section{Cylinder and Torus}
The Gross-Taylor formula is free of $\Omega$-factors whenever the number of holes 
in spacetime is equal to its Euler characteristic. The equation 
$m=2-2n$ does not have many solutions in nonnegative integers.

	\subsection{Cylinder}
	If $m=2$ and $n=0$, spacetime is a cylinder of area $t$ and the 
	microchiral partition function is
	
		\begin{equation}
		\label{eqn:MicrochiralCylinder}
				Z_N^d = \frac{z^d}{d!} \sum_{\alpha,\beta \in \Y^d} P_{\alpha\beta} 
				\sum_{\lambda \in \Y_N^d} \frac{(\dim \V^\lambda)^2}{d!} \hat{K}_{\alpha\beta}(\lambda) \hat{\Psi}_{\frac{t}{N}}(\lambda),
		\end{equation}		
		
	\noindent
	where $P_{\alpha\beta}(U_1,U_2) = p_\alpha(U_1) p_\alpha(U_2)$ is a product of trace
	invariants and $\hat{K}_{\alpha\beta}(\lambda) = \hat{K}_\alpha(\lambda) \hat{K}_\beta(\lambda)$
	is a product of central characters.
	For $1 \leq d \leq N$, the corresponding Gross-Taylor expectation is (two holes, zero handles)
	
		\begin{equation}
		\label{eqn:CylinderExpectation}
			\mathcal{E}_N^d(\alpha,\beta;0) =  \langle K_{\alpha\beta}\Psi_\frac{t}{N}\rangle.
		\end{equation}

	\subsubsection{Microchiral $1/N$ expansion}
	The cylinder expectation is 
	
		\begin{equation}
			\mathcal{E}_N^d(\alpha,\beta;0) = \sum_{r=0}^\infty \frac{(-t)^r}{r!N^r} 
			\langle K_{\alpha\beta} K^r \rangle,
		\end{equation}
		
	\noindent
	where a finite sum over $\lambda$ has been interchanged with an absolutely convergent series in $1/N$. 
	The Plancherel expectation
		
		\begin{equation}
			\langle K_{\alpha\beta}K^r \rangle = \langle K_\alpha K^r K_\beta \rangle = \langle K_\alpha | K^r | K_\beta \rangle		
		\end{equation}
		
	\noindent 
	counts walks on the Hurwitz-Cayley graph which begin at a point of the conjugacy class $K_\alpha$ and end at a point of 
	the conjugacy class $K_\beta$. We thus find that the cylinder expectation \eqref{eqn:CylinderExpectation} admits an absolutely 
	convergent series expansion in powers of $1/N$ whose coefficients enumerate walks between 
	specified conjugacy classes of $\group{S}^d$.
	
	The path count $\langle K_{\alpha\beta}K^r \rangle$ can be understood topologically 
	by inverting the monodromy correspondence \cite{CM}. More precisely, the normalized path count 
	$\frac{1}{d!} \langle K_{\alpha\beta}K^r\rangle$ counts orbits of the action of $\group{S}^d$
	on $r$-step walks $K_\alpha \to K_\beta$ in the Hurwitz-Cayley graph by simultaneous conjugation of steps and endpoints,
	each orbit being weighted by the reciprocal of the cardinality of the corresponding stabilizer.
	Length $r$ walks $K_\alpha \to K_\beta$ modulo conjugation are in bijection with 
	equivalence classes of pairs $(X,f)$ consisting of a 
	compact but not necessarily connected Riemann surface $X$ together with a degree
	$d$ holomorphic function $f \colon X \to Y$ to the Riemann sphere with a fixed
	branch locus $\{y_\alpha,y_\beta,y_1,\dots,y_r\}$ on $Y$ such that the 
	ramification profile of $f$ over the branch points $y_\alpha$ and $y_\beta$ is given by the Young 
	diagrams $\alpha$ and $\beta$, respectively, and the remaining branch points
	$y_1,\dots,y_r$ are simple.
	The normalized expectations $\frac{1}{d!} \langle K_{\alpha\beta}K^r\rangle$ 
	were called \emph{double Hurwitz numbers} in \cite{Okounkov}, owing to the fact that
	Hurwitz had considered these numbers in the case where one of $K_\alpha,K_\beta$ is
	the identity class; see \cite{Lando:ICM}. We will use the term ``double Hurwitz number''
	for the raw count $\langle K_{\alpha\beta} K^r\rangle$
	of $r$-step walks $K_\alpha \to K_\beta$ on the Hurwitz-Cayley graph.
					
	To conclude, inverting the usual flow of information
			
		\begin{equation}
		\begin{CD}
			\text{Surfaces} @>\text{Monodromy}>> \text{Permutations} @>\text{Fourier}>> \text{Characters},
		\end{CD}
		\end{equation}
		
	\noindent
	we find that the microchiral partition function of
	$YM_2$ on the cylinder is a generating function for double Hurwitz numbers.
			
		\begin{thm}
		\label{thm:MicrochiralCylinderExpansion}
			For any integers $1 \leq d \leq N$, the microchiral partition function $Z_N^d$
			of $YM_2$ on a cylinder of area $t$ admits the absolutely convergent series expansion 
			
				\begin{equation*}
					Z_N^d = \frac{z^d}{d!} \sum_{\alpha,\beta \in \Y^d} P_{\alpha\beta} 
					\sum_{r=0}^\infty \frac{(-t)^r}{r!N^r} \langle K_{\alpha\beta}  K^r\rangle,
				\end{equation*}
				
			\noindent
			where $\langle K_{\alpha\beta} K^r \rangle = \langle K_\alpha | K^r | K_\beta\rangle$ is a double Hurwitz number.
		\end{thm}
		
	Yang-Mills on the cylinder is particularly significant given its relation to the 
	heat kernel on the unitary group \cite{DGHK,KW,GrossMat,Zelditch}. 
	Double Hurwitz numbers are major objects of study in Hurwitz theory due to their deep connections
	with enumerative geometry and integrable systems \cite{BDKLM,GJV,KL,Okounkov}.
	However, the relationship between $YM_2$ on the cylinder and double Hurwitz numbers
	given by Theorem \ref{thm:MicrochiralCylinderExpansion} does not seem to be present in the existing
	literature on either subject. 
	
	Note that while Theorem \ref{thm:MicrochiralCylinderExpansion} pertains to $YM_2$ 
	with cylindrical spacetime, it is degree $d$ maps to the sphere which are counted by the $1/N$ expansion of $Z_N^d$. 
	Nevertheless, $1/N$ respects the boundaries of spacetime by remembering them as two distinguished points
	on the sphere, which we have taken to be $0$ and $\infty$, and counting maps to the sphere which may 
	have any ramification profile over these points. 
		
		\subsubsection{Chiral Gross-Taylor series}
		Theorem \ref{thm:MicrochiralCylinderExpansion} suggests that the large $N$ limit of $YM_2$ on the cylinder
		is double Hurwitz theory, the enumerative theory of maps from compact
		Riemann surfaces to the sphere with two branch points of arbitrary ramification, but it does not 
		directly imply such a statement in any quantitatively meaningful way.
		What we can do ``for free'' with Theorem \ref{thm:MicrochiralCylinderExpansion} 
		is formally set $N=\infty$. 
				
			\begin{definition}
			\label{eqn:ChiralSeries}
			The chiral Gross-Taylor series of the cylinder is
			the formal power series
				
					\begin{equation*}
						Z = 1 + \sum_{d=1}^\infty \frac{z^d}{d!} \sum_{\alpha,\beta \in \Y^d} P_{\alpha\beta}
						\sum_{r=0}^\infty \frac{(-t\hbar)^r}{r!} \langle K_{\alpha\beta} K^r\rangle,
					\end{equation*}
					
			\noindent
			where $z,t,\hbar,p_1,p_2,p_3,\dots$ are commuting indeterminates
			and 
		
			\begin{equation*}
				P_{\alpha\beta} = p_\alpha \otimes p_\beta = \left(\prod_{i=1}^{\ell(\alpha)} p_{\alpha_i} \right)
				\otimes \left(\prod_{i=1}^{\ell(\beta)} p_{\beta_i} \right).
			\end{equation*}
			\end{definition}
		
		It is important to understand that the chiral Gross-Taylor series is just a convenient
		way to organize the information in Theorem \ref{thm:MicrochiralCylinderExpansion}
		and does not have any additional meaning. However, 		
		$Z$ can be brought to a form which suggests
		what the large $N$ asymptotics of the chiral partition function $Z_N$ may actually be.
		The indeterminate $z$ is an exponential marker for the degree $d$ of the symmetric
		group, while $-t\hbar$ is an exponential marker for the length $r$ of a walk $K_\alpha \to K_\beta$ on 
		$\group{S}^d$ counted by $\langle K_{\alpha\beta} K^r \rangle$. From an enumerative
		perspective we could dispense
		with $\hbar$, but in order to anticipate asymptotics it is useful to keep it 
		as an infinitesimal form of $1/N$. The tensor $P_{\alpha\beta}$ is a multiplicative marker for the
		boundary conditions of a walk $K_\alpha \to K_\beta$ on $\group{S}^d$. 
		The Exponential Formula \cite{GJ:Book} gives us the formal logarithm of $Z$ as a generating 
		series enumerating ``connected'' walks on the Cayley graph.
		
			\begin{thm}
			\label{thm:TopologicalExpansionCylinder}
				The chiral Gross-Taylor series of the cylinder is given by $Z=e^F$, where
				
					\begin{equation*}
						F = \sum_{d=1}^\infty \frac{z^d}{d!} 
						\sum_{\alpha,\beta \in \Y^d} P_{\alpha\beta} 
						\sum_{r=0}^\infty \frac{(-t\hbar)^r}{r!} \langle K_{\alpha\beta} K^r \rangle_c
					\end{equation*}
					
				\noindent
				and the cumulant $\langle K_{\alpha\beta} K^r\rangle_c$
				is a connected double Hurwitz number.
			\end{thm}
						
		The combinatorial meaning of the connected double Hurwitz number $\langle K_{\alpha\beta}K^r \rangle$ is
		that it counts $r$-step walks $K_\alpha \to K_\beta$ on the Hurwitz-Cayley graph whose steps and endpoints together 
		generate a transitive subgroup of $\group{S}^d$. Topologically, $\frac{1}{d!}\langle K_{\alpha\beta} K^r  \rangle_c$ 
		enumerates degree $d$ covers
		as above, but with $X$ irreducible. By the Riemann-Hurwitz formula, $\langle K_{\alpha\beta} K^r \rangle_c$
		vanishes unless $r = 2g-2+\ell(\alpha)+\ell(\beta)$ with $g \geq 0$ the genus of $X$. Setting
		
			\begin{equation}
				H_g(\alpha,\beta) =  \langle K_{\alpha\beta}  K^{2g-2+\ell(\alpha)+\ell(\beta)} \rangle_c,
			\end{equation}

		\noindent
		the logarithm of the chiral Gross-Taylor series becomes 
		
			\begin{equation}
				F = \sum_{d=1}^\infty \frac{z^d}{d!} \sum_{\alpha,\beta \in \Y^d} P_{\alpha\beta}
				\sum_{g=0}^\infty \frac{(-t\hbar)^{2g-2+\ell(\alpha)+\ell(\beta)}}{(2g-2+\ell(\alpha)+\ell(\beta))!} H_g(\alpha,\beta).
			\end{equation}
			
		\noindent
		This formal power series is the main object of study in \cite{GJV,Okounkov}.
		Since it is a formal series, we are free to sum first over genus and then over degree to get a ``genus expansion.''

		\begin{thm}
		\label{thm:StableCylinderExpansion}
			The chiral Gross-Taylor series $Z$ of $YM_2$ on a cylinder admits 
			the topological expansion 
			
				\begin{equation*}
					Z = \exp \sum_{g=0}^\infty \hbar^{2g-2} F_g,
				\end{equation*}
				
			\noindent
			where $F_g$ is given by
			
				\begin{equation*}
					F_g = t^{2g-2} \sum_{d=1}^\infty \frac{z^d}{d!} \sum_{\alpha,\beta \in \Y^d}
					P_{\alpha\beta} \frac{(-t\hbar)^{\ell(\alpha)+\ell(\beta)}}{(2g-2+\ell(\alpha)+\ell(\beta))!} H_g(\alpha,\beta),
				\end{equation*} 
				
			\noindent
			a generating series for connected double Hurwitz numbers of genus $g$.
		\end{thm}

		\subsubsection{Large $N$}
		Theorem \ref{thm:StableCylinderExpansion} is a formal statement, but it offers
		a prediction on the $N \to \infty$ behavior of the chiral partition 
		function 
		
			\begin{equation}
				Z_N(U_1,U_2) = 1 + \sum_{\lambda \in \Y_N} s_\lambda(U_1)s_\lambda(U_2) e^{-\frac{t}{2N}\hat{C}(\lambda)}
			\end{equation}
		
		\noindent
		of $YM_2$ on a cylinder of area $t$. 
		Remembering that $\hbar = 1/N$ with $N=\infty$, it 
		is natural to conjecture on the basis of Theorem \ref{thm:StableCylinderExpansion} that
				
			\begin{equation}
			\label{eqn:AsymptoticExpansion}
				\log Z_N(U_1,U_2) \sim \sum_{g=0}^\infty N^{2-2g} F_g(U_1,U_2)
			\end{equation}
			
		\noindent
		as $N \to \infty$, where
		
			\begin{equation}
			\label{eqn:CylinderFreeEnergy}
				F_g(U_1,U_2) = t^{2g-2}\sum_{d=1}^\infty \frac{e^{-\frac{td}{2}}}{d!}
				\sum_{\alpha,\beta \in \Y_N^d} \frac{p_\alpha(U_1)}{N^{\ell(\alpha)}} \frac{p_\beta(U_2)}{N^{\ell(\beta)}}
				 \frac{(-t)^{\ell(\alpha)+\ell(\beta)}}{(2g-2+\ell(\alpha)+\ell(\beta))!} H_g(\alpha,\beta).
			\end{equation}
		
		\noindent
		This putative $N \to \infty$ asymptotic expansion of $\log Z_N$ is 
		inspired by Theorem \ref{thm:StableCylinderExpansion}, not implied by it.
		Indeed, one cannot even claim that the right hand side of \eqref{eqn:AsymptoticExpansion}
		is an asymptotic series without first demonstrating that coefficients $F_g$ defined 
		by \eqref{eqn:CylinderFreeEnergy} converge for $t$ in some non-trivial and $g$-independent set 
		$\mathcal{T}$ of positive number (this convergence set is related to the behavior of Brownian motion on 
		$\group{U}_N$ as $N \to \infty$). We will discuss this 
		further below for $YM_2$ with spherical spacetime, where the relationship between
		the chiral Gross-Taylor series and a true large $N$ approximation is even clearer.
											
	\subsection{Torus}
	If $m=0$ and $n=1$, spacetime is a torus of area $t$ and the 
	microchiral partition function is 
	
		\begin{equation}
		\label{eqn:MicrochiralTorus}
				Z_N^d =  \frac{z^d}{d!}
				\sum_{\lambda \in \Y_N^d}\frac{(\dim \V^\lambda)^2}{d!}
				\hat{H}(\lambda) \hat{\Psi}_\frac{t}{N}(\lambda).
			\end{equation}

	\noindent
	For $1 \leq d \leq N$, the corresponding
	 Gross-Taylor expectation is (no holes, one handle)
		
		\begin{equation}
		\label{eqn:TorusExpectation}
			\mathcal{E}_N^d(1) =  \langle H\Psi_\frac{t}{N}\rangle.
		\end{equation}
		
	\noindent
	Since a torus has no boundaries there are no Young diagrams in the 
	expectation \eqref{eqn:TorusExpectation}.

	\subsubsection{Microchiral $1/N$ expansion}
	The torus expectation is 
	
		\begin{equation}
			\mathcal{E}_N^d(1) = \sum_{r=0}^\infty \frac{(-t)^r}{r!N^r} 
			\langle HK^r \rangle,
		\end{equation}
		
	\noindent
	where a finite sum over $\lambda$ has been interchanged with an absolutely convergent series in $1/N$. 
	The Plancherel expectation $\langle HK^r\rangle$ counts the number of solutions of the equation 
		
		\begin{equation}
		\label{eqn:BrokenStick}
			\rho\sigma\tau_1 \dots \tau_r = \rho\sigma		
		\end{equation}
		
	\noindent 
	in the symmetric group $\group{S}^d$ such that $\rho,\sigma$ are arbitrary and the $\tau_i$'s
	are transpositions.

	Solutions of the equation \eqref{eqn:BrokenStick} have a well-known topological 
	interpretation which we will discuss momentarily, but first we explain how they can 
	be interpreted graph-theoretically which does not seem to be well-known. Invoking
	the Jucys-Murphy correspondence, every two-factor product $\rho\sigma$ corresponds 
	uniquely to a concatenation of two strictly monotone walks $\rho$ and $\sigma$ on the Hurwitz-Cayley graph
	such that the initial point of $\rho$ is $\iota$ and the initial point of $\sigma$ is the terminal point of $\rho$.
	Therefore, the Plancherel expectation $\langle HK^r \rangle$ is the number of walks on 
	on the Hurwitz-Cayley graph which begin at the identity $\iota$ and consist of a concatenation of two strictly monotone walks
	$\rho$ and $\sigma$, each of any length from $0$ to $d-1$, followed by loop with steps $\tau_1,\dots,\tau_r$ 
	based at $\rho\sigma$.
		
	The topological interpretation of $\langle HK^r\rangle$ is classical: after normalizing
	by $\frac{1}{d!}$ this expectation counts isomorphism classes of degree $d$ maps $f \colon X \to Y$ from a compact 
	but not necessarily connected Riemann surface $X$ to the torus $Y$ which have 
	$r$ simple branch points $y_1,\dots,y_r$ on $Y$ and no other branching \cite{Zagier}.
	For this reason $\frac{1}{d!}\langle HK^r \rangle$ is called a \emph{simple Hurwitz number} with 
	torus target, but we will use this term for the raw broken-stick loop count $\langle HK^r \rangle$.
		
	\begin{thm}
		\label{thm:MicrochiralTorusExpansion}
			For any integers $1 \leq d \leq N$, the microchiral partition function $Z_N^d$
			of $YM_2$ on a torus of area $t$ admits the absolutely convergent series expansion 
			
				\begin{equation*}
					Z_N^d = \frac{z^d}{d!}
					\sum_{r=0}^\infty \frac{(-t)^r}{r!N^r} \langle HK^r\rangle,
				\end{equation*}
				
			\noindent
			where $\langle HK^r\rangle$ is a simple Hurwitz number with torus target.
		\end{thm}
		
		\subsubsection{Gross-Taylor series}
		
			\begin{definition}
				The chiral Gross-Taylor series of the torus is the trivariate formal power series
				
					\begin{equation*}
						Z = \sum_{d=0}^\infty \frac{z^d}{d!} \sum_{r=0}^\infty \frac{(-t\hbar)^r}{r!} 
						\langle HK^r \rangle.
					\end{equation*}
			\end{definition}

		Once again, $z$ is an exponential marker for degree and $-t\hbar$ is an exponential 
		marker for the cardinality of the branch locus of a map with simple 
		branch points at prescribed locations, or equivlanelty the length of a loop in the Hurwitz-Cayley graph $\group{S}^d$
		attached to a broken stick. The Exponential Formula gives the formal logarithm of $Z$
		as the solution to a connected enumeration problem.
		
			\begin{thm}
				We have $Z= e^F$, where
				
					\begin{equation*}
						F = \sum_{d=1}^\infty \frac{z^d}{d!} \sum_{r=0}^\infty \frac{(-t\hbar)^r}{r!} 
						\langle HK^r\rangle_c,
					\end{equation*}
					
				\noindent
				and the cumulant $\langle HK^r\rangle_c$ is a connected
				simple Hurwitz number over the torus.
			\end{thm}
	
		By the Riemann-Hurwitz formula, the cumulant $\langle HK^r\rangle_c$ vanishes
		unless $r=2g-2$ with $g \geq 1$ the genus of $X$, so we set
		
			\begin{equation}
				H_g^d = \langle HK^{2g-2} \rangle_c.
			\end{equation}
			
		\noindent
		Then, we have the following genus expansion of the chiral Gross-Taylor series.
			
			\begin{thm}
			\label{thm:StableTorusExpansion}
			The chiral Gross-Taylor series $Z$ of $YM_2$ on a torus admits 
			the topological expansion 
			
				\begin{equation*}
					Z = \exp \sum_{g=1}^\infty \hbar^{2g-2} F_g,
				\end{equation*}
				
			\noindent
			where $F_g$ with
			
				\begin{equation*}
					F_g= \frac{t^{2g-2}}{(2g-2)!}\sum_{d=1}^\infty \frac{z^d}{d!} H_g^d.
				\end{equation*} 
			\end{thm}
		
		The univariate formal power series $\tilde{F}_g = \sum_{d=1}^\infty \frac{z^d}{d!} H_g^d$
		is the main object of study in \cite{Dijkgraaf}, where it is shown that for $g \geq 2$ the 
		series $\tilde{F}_g$ is a quasimodular form of weight $6g-6$, i.e. a polynomial in the 
		Eisenstein series $E_2,E_4,E_6$ of degree $6g-6$ with respect to the grading
		$\deg E_k = k$. We refer to \cite{Dijkgraaf} for a full discussion of this remarkable
		fact, which gives a way to understand the convergence properties of the 
		genus $g$ free energies $F_g$ in terms of properties of quasimodular forms. 
						
\section{Sphere and Disc}
The Gross-Taylor formula contains a positive power of $\Omega$ whenever the number of holes 
in spacetime is less than its Euler characteristic. The inequality 
$m<2-2n$ does not have many solutions
in nonnegative integers.

	\subsection{Sphere}
	If $m=0$ and $n=0$, spacetime is a sphere of area $t$ and the $YM_2$
	microchiral partition function is 
	
		\begin{equation}
		\label{eqn:MicrochiralSphere}
				Z_N^d = \frac{(zN^2)^d}{d!}
				\sum_{\lambda \in \Y_N^d} \frac{(\dim \V^\lambda)^2}{d!} \hat{\Psi}_{\frac{t}{N}}(\lambda)
				\hat{\Omega}_\frac{1}{N}^2(\lambda).
		\end{equation}		
		
	\noindent
	The corresponding Gross-Taylor expectation is (no holes, no handles) 
	
		\begin{equation}
		\label{eqn:SphereExpectation}
			\mathcal{E}_N^d = \langle \Psi_\frac{t}{N} \Omega_\frac{1}{N}^2\rangle.
		\end{equation}
				
		\subsubsection{Microchiral partition function}
		Both sources of $1/N$ are now in play. Write the Gross-Taylor expectation 
		\eqref{eqn:SphereExpectation} as $\langle \Omega_\frac{1}{N} \Psi_\frac{t}{N} \Omega_\frac{1}{N} \rangle$
		and recall from Section \ref{sec:HurwitzCayley} that the exponential distance element in 
		$\C\group{S}^d$ is 
		
			\begin{equation}
				|\Omega_q\rangle = \sum_{a=0}^{d-1} q^a |L_a\rangle
			\end{equation}
			
		\noindent
		with $|L_a\rangle$ the $a$th level of the Hurwitz-Cayley graph $\group{S}^d$. Then 
		
			\begin{equation}
				\langle \Omega_\frac{1}{N} \Psi_\frac{t}{N} \Omega_\frac{1}{N} \rangle = 
				\sum_{a,b=0}^{d-1} \frac{1}{N^{a+b}} \sum_{r=0}^\infty \frac{(-t)^r}{r!N^r}
				\langle L_a K^r L_b \rangle,
			\end{equation}
			
		\noindent
		where the Plancherel expectation 
		
			\begin{equation}
				\langle L_a K^r L_b \rangle = \langle L_a | K^r | L_b \rangle 
			\end{equation}
			
		\noindent
		counts $r$-step walks $L_a \to L_b$ in the Hurwitz-Cayley graph.
		This is a sum of double Hurwitz numbers with spherical target, 
				
			\begin{equation}
				\langle L_a K^r L_b \rangle = \sum_{\substack{\alpha,\beta \in \Y^d \\
				d-\ell(\alpha)=a,\ d-\ell(\beta)=b}} \langle K_\alpha K^r K_\beta \rangle,
			\end{equation}
			
		\noindent
		and instead of depending on a pair of partitions $\alpha,\beta \in \Y^d$ marking 
		conjugacy classes in $\group{S}^d$ it depends on a pair of integers $a,b \in \{0,1,\dots,d-1\}$
		marking levels in $\group{S}^d$. We therefore call $\langle L_a K^r L_b \rangle$
		a \emph{coarse double Hurwitz number} over the sphere.

			\begin{thm}
			\label{thm:MicrochiralSphereExpansion}
				For any integers $1 \leq d \leq N$, the microchiral partition function $Z_N^d$
				of $YM_2$ on a sphere of area $t$ admits the absolutely convergent series expansion 
			
					\begin{equation*}
						Z_N^d = \frac{(zN^2)^d}{d!} \sum_{a,b=0}^{d-1} \frac{1}{N^{a+b}}
						\sum_{r=0}^\infty \frac{(-t)^r}{r!N^r} \langle L_a K^r L_b \rangle,
					\end{equation*}
				
				\noindent
				where $\langle L_a K^r L_b \rangle$ is a coarse double Hurwitz number.
			\end{thm}
			
		Considerable effort has gone into trying to guess the form of the $1/N$ expansion
		of $YM_2$ on the sphere \cite{CMR1,CMR2,Taylor}. While 
		Theorem \ref{thm:MicrochiralSphereExpansion} gives a mathematically simple answer, 
		the physical meaning of this answer is non-obvious.
		On one hand, the case of a spherical 
		spacetime is like the case of a torus spacetime in that the $1/N$ 
		expansion actually enumerates worldsheet maps 
		to spacetime. On the other hand, the spherical spacetime seems
		to conceive of itself as a cylinder with two nonexistent boundaries which manifest as two
		distinguished points over which the ramification type of the map is arbitrary;
		the $1/N$ expansion sees only the number of cycles 
		in the ramification profile of these points.

		\subsubsection{Chiral Gross-Taylor series}
		Since the chiral Gross-Taylor series is by definition a formal series, we may 
		enrich it with additional parameters which are useful enumerative markers not present in the
		quantitative Gross-Taylor formula.
		
			\begin{definition}
			\label{def:ChiralGrossTaylorSphere}
				The chiral Gross-Taylor series of the sphere is defined by
		
				\begin{equation*}
					Z = 1 + \sum_{d=1}^\infty \frac{(z\hbar^{-2})^d}{d!}
					\sum_{a,b=0}^{d-1} (u\hbar)^a (v\hbar)^b \sum_{r=0}^\infty \frac{(-t\hbar)^r}{r!} 
					\langle L_a K^r L_b \rangle,
				\end{equation*}
				
			\noindent
			where $z,\hbar,t,u,v$ are formal variables.
			\end{definition}

		As in the case of the chiral Gross-Taylor series of the cylinder, 
		the indeterminate $\hbar$ seems to be redundant, but ultimately it serves
		a purpose. The indeterminate $z$ is once again
		an exponential marker for the degree of $\group{S}^d$, and $t$ is again 
		an exponential marker for the number of steps in a length $r$ walk $L_a \to L_b$.
		The indeterminates $u$ and $v$ are ordinary markers for the boundary
		conditions $a,b \in \{0,\dots,d-1\}$. 
		The Jucys-Murphy correspondence gives a second interpretation of 
		the expectation $\langle L_a K^r L_b \rangle$ as counting loops of length
		$a+r+b$ on the Hurwitz-Cayley graph based at $\iota$ which consist of $a$ strictly monotone increasing steps, followed
		by $r$ unrestricted steps, followed by $b$ strictly monotone decreasing steps.
		the monotone walks which make up the loop carry ordinary weight because their steps cannot be shuffled,
		whereas unrestricted walks carry exponential weight because their steps can be shuffled.
						
			\begin{thm}
				We have $Z=e^F$, where 
				
					\begin{equation*}
						F = \sum_{d=1}^\infty \frac{(z\hbar^{-2})^d}{d!}
						\sum_{a,b=0}^{d-1} (u\hbar)^a (v\hbar)^b \sum_{r=0}^\infty \frac{(-t\hbar)^r}{r!} 
						\langle L_a K^r L_b \rangle_c
					\end{equation*}
					
				\noindent
				and the cumulant $\langle L_a K^r L_b \rangle_c$ is the number
				of $r$-step walks $L_a \to L_b$ on the Hurwitz-Cayley graph whose steps and 
				endpoints generate a transitive subgroup of $\group{S}^d$. 			
			\end{thm}
		
		The cumulant $\langle L_a K^r L_b \rangle_c$ counts connected branched
		covers with $d-a$ and $d-b$ cycles in the ramification over two branch 
		points and the remaining $r$ branch points simple. From the
		Riemann-Hurwitz formula, we see that the cumulant $\langle L_a K^r L_b \rangle_c$
		is zero unless 
		
			\begin{equation}
			\label{eqn:CoarseRiemannHurwitz}
				r = 2g-2 + (d-a) + (d-b)
			\end{equation}
			
		\noindent
		for a nonnegative integer $g$, the genus of the cover.
		Therefore, it is natural to parameterize the connected coarse double Hurwitz numbers
		by genus, writing 
		
			\begin{equation}
			\label{eqn:ConnectedCoarseDouble}
				H_g(a,b) = \langle L_a K^{2g-2+2d-a-b} L_b \rangle.
			\end{equation}
			
		\noindent
		In the genus parameterization, the net power of $\hbar$ in the coefficient 
		of $z^d$ in $F$ is
		
			\begin{equation}
				\hbar^{-2d} \hbar^a \hbar^b \hbar^{2g-2+2d-a-b} = \hbar^{2g-2}.
			\end{equation}
		
		\noindent
		This gives a genus expansion of the chiral Gross-Taylor free energy 
		of $YM_2$ on the sphere in the physical specialization $u=v=1$ corresponding
		to the expectation value $\langle \Omega_\frac{1}{N} \Psi_\frac{t}{N} \Omega_\frac{1}{N} \rangle$.
		
			\begin{thm}
			\label{thm:ChiralFreeEnergySphere}
				We have 
				
					\begin{equation*}
						F \big{|}_{u=v=1} = \sum_{g=0}^\infty \hbar^{2g-2} F_g,
					\end{equation*}
					
				\noindent
				where 
				
					\begin{equation*}
						F_g = t^{2g-2}\sum_{d=1}^\infty \frac{z^d}{d!} \sum_{a,b=0}^{d-1} t^{2d-a-b}
							\frac{H_g(a,b)}{(2g-2+2d - a - b)!}.
					\end{equation*}
			\end{thm}
			
		\subsubsection{Large $N$}
		As in the case of the cylinder, Theorem \ref{thm:ChiralFreeEnergySphere} predicts an $N \to \infty$ approximation 
		of the chiral partition function $Z_N$ of $YM_2$ on a sphere of area $t>0$. The predicted
		asymptotics are
		
			\begin{equation}
			\label{eqn:SpherePrediction}
				\log Z_N \sim \sum_{g=0}^\infty N^{2-2g} F_g,
			\end{equation}
			
		\noindent
		where 
		
			\begin{equation}
			\label{eqn:GenusSpecificSphere}
				F_g = t^{2g-2}\sum_{d=1}^\infty \frac{e^{-\frac{td}{2}}}{d!} \sum_{a,b=0}^{d-1} t^{2d-a-b}
							\frac{H_g(a,b)}{(2g-2+2d - a - b)!}
			\end{equation}
			
		\noindent
		is a positive series. In order for the prediction \ref{eqn:SpherePrediction} to 
		even make sense, it is necessary that there exists a set $\mathcal{T} \subseteq \R_{\geq 0}$ the series $F_g=F_g(t)$ given by 
		\eqref{eqn:GenusSpecificSphere} converges for every $g \in \N_0$. 
		The formal power series
		
			\begin{equation}
				\tilde{F}_g = \sum_{d=1}^\infty \frac{z^d}{d!} \sum_{a,b=0}^{d-1} t^{2d-a-b}
							\frac{H_g(a,b)}{(2g-2+2d - a - b)!}
			\end{equation}
			
		\noindent
		with $z$ independent of $t$
		does indeed have radius of convergence bounded below by a positive constant 
		independent of $g$, see \cite{GGN5} and \cite{BCCG}. Setting $z=e^{-t/2}$ with $t >0$ means that 
		the convergence set $\mathcal{T} \subseteq \R_{\geq 0}$ consists of two disjoint
		intervals, as predicted in \cite{Taylor}. This is related to the large $N$ behavior of Brownian 
		loops in $\group{U}_N$.

	\subsection{Disc}
	In addition to the sphere, there is one and only one other spacetime for which the Gross-Taylor formula
	contains a positive power of $\Omega$.
	If $m=1$ and $n=0$, spacetime is a disc of area $t$ and the $YM_2$
	microchiral partition function is 
	
		\begin{equation}
		\label{eqn:MicrochiralDisc}
				Z_N^d = \frac{(zN)^d}{d!} \sum_{\alpha \in \Y^d} p_\alpha
				\sum_{\lambda \in \Y_N^d} \frac{(\dim \V^\lambda)^2}{d!} \hat{K}_\alpha(\lambda) \hat{\Psi}_{\frac{t}{N}}(\lambda)
				\hat{\Omega}_\frac{1}{N}(\lambda).
		\end{equation}
		
	\noindent
	The corresponding Gross-Taylor expectation is (one hole, no handles)
	
		\begin{equation}
		\label{eqn:DiscExpectation}
			\mathcal{E}_N^d(\alpha;0) = \langle K_\alpha \Psi_\frac{t}{N} \Omega_\frac{1}{N} \rangle.
		\end{equation}
		
	\noindent
	The $1/N$ expansion of \eqref{eqn:DiscExpectation} has coefficients given by the 
	Plancherel expecations
	$\langle K_\alpha K^r L_b \rangle$, which are ``half-coarse'' double Hurwitz numbers
	counting $r$-step walks $K_\alpha \to L_b$ from a specified conjugacy class to a specified 
	level in the Hurwitz-Cayley graph, or equivalently branched covers of the sphere
	whose profile over $0$ is $\alpha$ while the ramification over $\infty$ has $d-b$ cycles.
	Being almost identical to what has come before, the details are left to the interested reader.
	
\section{Three-Holed Sphere and One-Holed Torus}
The Gross-Taylor formula contains a negative power of $\Omega$ whenever the number of holes 
in spacetime is greater than its Euler characteristic. The inequality 
$m>2-2n$ has infinitely many solutions
in nonnegative integers. We now examine 
the two extremal solutions in which the number of holes in spacetime
exceeds the Euler characteristic of spacetime by exactly one.

	\subsection{Three-holed sphere}
	If $m=3$ and $n=0$, spacetime is a size $t$ pair of pants and the $YM_2$
	microchiral partition function is 
	
		\begin{equation}
		\label{eqn:MicrochiralSphere}
				Z_N^d = \frac{(zN^{-1})^d}{d!} \sum_{\alpha,\beta,\gamma \in \Y^d} P_{\alpha\beta\gamma}
				\sum_{\lambda \in \Y_N^d} \frac{(\dim \V^\lambda)^2}{d!} \hat{K}_{\alpha\beta\gamma}(\lambda)
				\hat{\Psi}_{\frac{t}{N}}(\lambda)\hat{\Omega}_\frac{1}{N}^{-1}(\lambda),
		\end{equation}
		
	\noindent
	where $P_{\alpha\beta\gamma}(U_1,U_2,U_3) = p_\alpha(U_1) p_\beta(U_2) p_\gamma(U_3)$ is a 
	product of three trace invariants and $\hat{K}_{\alpha\beta\gamma}(\lambda) =
	\hat{K}_\alpha(\lambda)\hat{K}_\beta(\lambda)\hat{K}_\gamma(\lambda)$
	is a product of three central characters. The corresponding Gross-Taylor expectation for $1 \leq d \leq N$ is
	(three holes, no handles)
	
		\begin{equation}
		\label{eqn:HoleySphereExpectation}
			\mathcal{E}_N^d(\alpha,\beta,\gamma;0) = \langle K_{\alpha\beta\gamma} \Psi_\frac{t}{N} \Omega_\frac{1}{N}^{-1} \rangle.
		\end{equation}

		\subsubsection{Microchiral partition function}
		In order to obtain the $1/N$ expansion of the Gross-Taylor expectation \eqref{eqn:HoleySphereExpectation},
		we must use the result from Section \ref{sec:HurwitzCayley} that the exponential distance operator
		$\Omega_q$ on the Hurwitz-Cayley graph $\group{S}^d$ is invertible for $|q| < \frac{1}{d-1}$ with
		inverse given by the absolutely convergent series
		
			\begin{equation}
				\Omega_q^{-1} = \sum_{r=0}^\infty (-q)^r M_r,
			\end{equation}
			
		\noindent
		where $M_r=h_r(J_1,\dots,J_d)$ is the complete symmetric polynomial of degree $r$
		in the Jucys-Murphy operators. Thus for $1 \leq d \leq N$ we have the absolutely convergent
		$1/N$ expansion 
		
			\begin{equation}
				\mathcal{E}_N^d(\alpha,\beta,\gamma;0) = \sum_{r,s=0}^\infty \frac{(-t)^r (-1)^s}{r!N^{r+s}} 
				\langle K_{\alpha\beta\gamma} K^r M_s \rangle,
			\end{equation}

		\noindent
		where the Plancherel expectation $\langle K_{\alpha\beta\gamma} K^r M_s \rangle$
		counts $(3+r+s)$-tuples $(\pi,\rho,\sigma,\tau_1,\dots,\tau_r,\mu_1,\dots,\mu_s)$ from the 
		symmetric group $\group{S}^d$ such that
		
			\begin{equation}
				\pi\rho\sigma\tau_1 \dots \tau_r \mu_1\dots \mu_s=\iota,
			\end{equation}
			
		\noindent
		with 
		
			\begin{equation}
				\langle \pi | K_\alpha \rangle = \langle \rho | K_\beta\rangle = \langle \sigma | K_\gamma\rangle =1 
			\end{equation}
			
		\noindent
		and $\tau_1,\dots,\tau_r,\mu_1,\dots,\mu_s \in K$ transpositions, such that moreover 
		
			\begin{equation}
				\mu_1 \dots \mu_s = (i_1 j_1) \dots (i_s\ j_s), \quad i_k < j_k,\ j_1 \leq \dots \leq j_s.
			\end{equation}
			
		The expectation $\langle K_{\alpha\beta\gamma} K^r M_s \rangle$ is a \emph{mixed triple Hurwitz number}.
		To interpret it graphically as counting walks in the Hurwitz-Cayley graph, we first 
		decompose 
		
			\begin{equation}
				K_{\alpha\beta\gamma} = K_\alpha K_\beta K_\gamma = \sum_{\eta \in \Y^d}
				c_{\beta\gamma}^\eta K_\alpha K_\eta,
			\end{equation}
			
		\noindent
		where $c_{\beta\gamma}^\eta$ are the \emph{connection coefficients} of the class algebra,
		i.e. the structure constants of the center of $\C\group{S}^d$ relative the class basis, which 
		are nonnegative integers. This gives the mixed triple Hurwitz number as the nonnegative 
		integral linear combination 
		
			\begin{equation}
				\langle K_{\alpha\beta\gamma} K^r M_s \rangle =  \sum_{\eta \in \Y^d}
				c_{\beta\gamma}^\eta \langle K_\alpha K^r M_s K_\eta\rangle,
			\end{equation}
			
		\noindent
		where now each term $\langle K_\alpha K^r M_s K_\eta\rangle$ is a 
		\emph{mixed double Hurwitz number} \cite{GGN4} which counts length $r+s$
		walks $K_\alpha \to K_\eta$ in the Hurwitz-Cayley graph whose first $r$
		streps are unrestricted and whose last $s$ steps are weakly monotone
		in the sense of Section \ref{sec:HurwitzCayley}. In this very simple-minded
		sense, the connection coefficients of the class algebra can be viewed as
		giving an analogue of the Wick formula which allows to write an arbitrary Plancherel
		expectation $\langle K_{\alpha_1 \dots \alpha_m} \dots \rangle$ as a
		nonnegative integral linear combination of Plancherel expectations
		$\langle K_{\alpha\beta} \dots \rangle$.
		
		In the case $s=0$, the mixed triple Hurwitz number $\langle K_{\alpha\beta\gamma} K^r M_s\rangle$
		becomes the purely classical triple Hurwitz number $\langle K_{\alpha\beta\gamma}K^r\rangle$
		whose topological meaning is the same as that of the classical double Hurwitz numbers, 
		except that the covers $f \colon X \to Y$ of the sphere being counted have 
		three points with arbitrary ramification profile in their branch locus. In the 
		general case, one can say that the mixed triple Hurwitz number 
		$\langle K_{\alpha\beta\gamma} K^r M_s\rangle$ counts some subset of 
		the covers counted by the classical triple Hurwitz number $\langle K_{\alpha\beta\gamma} K^{r+s}\rangle$.
		A reasonable objection to this is that the subset of 
		covers counted depends on a labeling of the points in the fibre over the
		unbranched basepoint used to perform the monodromy construction, and 
		so is not geometrically defined. A reasonable response
		is that the \emph{cardinality} of the restricted subset of 
		covers does \emph{not} depend on the choice of a labeling because of 
		the centrality of symmetric polynomial functions of Jucys-Murphy elements.

			\begin{thm}
			\label{thm:MicrochiralThreeHoledSphereExpansion}
				For any integers $1 \leq d \leq N$, the microchiral partition function $Z_N^d$
				of $YM_2$ on a sphere of area $t$ admits the absolutely convergent series expansion 
			
					\begin{equation*}
						Z_N^d = \frac{(zN^{-1})^d}{d!} \sum_{\alpha,\beta,\gamma \in \Y^d} P_{\alpha\beta\gamma}
						\sum_{r,s=0}^\infty \frac{(-t)^r (-1)^s}{r!N^{r+s}} \langle K_{\alpha\beta\gamma} K^r M_s \rangle,
					\end{equation*}
				
				\noindent
				where $\langle K_{\alpha\beta\gamma} K^r M_s \rangle$ is a mixed triple Hurwitz number.
			\end{thm}
			
		\subsubsection{Chiral Gross-Taylor series}
		
			\begin{definition}
				The chiral Gross-Taylor series of the three-holed sphere is 
				
					\begin{equation*}
						Z= \sum_{d=0}^\infty \frac{(z\hbar)^d}{d!} \sum_{\alpha,\beta,\gamma \in \Y^d} 
						P_{\alpha\beta\gamma} \sum_{r,s=0}^\infty \frac{ (-t\hbar)^r}{r!} (-u\hbar)^s
						\langle K_{\alpha\beta\gamma} K^r M_s \rangle,
					\end{equation*}
					
				\noindent
				where $z,t,\hbar,p_1,p_2,p_3,\dots$ are indeterminates and 
				
					\begin{equation*}
						P_{\alpha\beta\gamma} = \left( \prod_{i=1}^{\ell(\alpha)} p_{\alpha_i}\right)
						\otimes \left( \prod_{i=1}^{\ell(\beta)} p_{\beta_i}\right) \otimes \left( \prod_{i=1}^{\ell(\gamma)} p_{\gamma}\right).
					\end{equation*}
			\end{definition}
 
 			\begin{thm}
				We have $Z= e^F$, where
				
					\begin{equation*}
						F = \sum_{d=0}^\infty \frac{(z\hbar^{-1})^d}{d!} \sum_{\alpha,\beta,\gamma \in \Y^d} 
						P_{\alpha\beta\gamma} \sum_{r,s=0}^\infty \frac{ (-t\hbar)^r}{r!} (-u\hbar)^s
						\langle K_{\alpha\beta\gamma} K^r M_s \rangle
					\end{equation*}
					
				\noindent
				and the cumulant $\langle K_{\alpha\beta\gamma} K^r M_s \rangle_c$
				is a connected mixed triple Hurwitz number.			
			\end{thm}
			
			From the Riemann-Hurwitz formula, $\langle K_{\alpha\beta\gamma} K^r M_s \rangle_c$
			is zero unless 
			
				\begin{equation}
					r+s=2g-2+d+\ell(\alpha)+\ell(\beta)+\ell(\gamma),
				\end{equation}
				
			\noindent
			so that the net power of $\hbar$ at degree $d$ in $F$ is
			
				\begin{equation}
					\hbar^d \hbar^{2g-2-d+\ell(\alpha)+\ell(\beta)+\ell(\gamma)} = \hbar^{2g-2} \hbar^{\ell(\alpha)+\ell(\beta)+\ell(\gamma)},
				\end{equation}
				
			\noindent
			and we obtain the topological expansion of the 
			chiral Gross-Taylor series of the three-holed sphere.

			\begin{thm}
			\label{thm:StableHoleySphereExpansion}
			The chiral Gross-Taylor series $Z$ of $YM_2$ on a three-holed sphere admits 
			the topological expansion 
			
				\begin{equation*}
					Z = \exp \sum_{g=0}^\infty \hbar^{2g-2} F_g,
				\end{equation*}
				
			\noindent
			where $F_g$ is given by
			
				\begin{equation*}
					F_g = \sum_{d=1}^\infty \frac{z^d}{d!} \sum_{\alpha,\beta,\gamma \in \Y^d}
					P_{\alpha\beta\gamma} \hbar^{\ell(\alpha)+\ell(\beta)+\ell(\gamma)} 
					\sum_{r+s=2g-2-d+\ell(\alpha)+\ell(\beta)+\ell(\gamma)} \frac{(-t)^r}{r!} (-u)^s
					\langle K_{\alpha\beta\gamma} K^r M_s \rangle_c.
				\end{equation*}

			\noindent
			a generating series for connected double Hurwitz numbers of genus $g$.
			\end{thm}

			\noindent
			Note that the internal sum of $F_g$ is finite, and its extreme terms correspond to connected classical and monotone 
			triple Hurwitz numbers: the $s=0$ term is 
			
				\begin{equation}
					\frac{(-t)^{\ell(\alpha)+\ell(\beta)+\ell(\gamma)}}{(2g-2-d+\ell(\alpha)+\ell(\beta)+\ell(\gamma))!}
					H_g(\alpha,\beta,\gamma)
				\end{equation}
				
			\noindent
			with 
			
				\begin{equation}
					H_g(\alpha,\beta,\gamma) = \langle K_{\alpha\beta\gamma}K^{2g-2-d+\ell(\alpha)+\ell(\beta)+\ell(\gamma)}\rangle_c
				\end{equation}
				
			\noindent
			a connected classical triple Hurwitz number, and the $r=0$ term is

				\begin{equation}
					(-u)^{\ell(\alpha)+\ell(\beta)+\ell(\gamma)}
					H_g^{\leq}(\alpha,\beta,\gamma)
				\end{equation}
				
			\noindent
			with 
			
				\begin{equation}
					H_g^{\leq}(\alpha,\beta,\gamma) = \langle K_{\alpha\beta\gamma}M_{2g-2-d+\ell(\alpha)+\ell(\beta)+\ell(\gamma)}\rangle_c
				\end{equation}
				
			\noindent
			a connected monotone triple Hurwitz number.

	\subsection{One-holed torus}
	If $m=1$ and $n=1$, spacetime is a one-holed torus of area $t$ and the 
	$YM_2$ microchiral partition function is 
	
		\begin{equation}
		\label{eqn:MicrochiralOneHoledTorus}
			Z_N^d = \frac{(zN^{-1})^d}{d!} \sum_{\alpha \in \Y^d} p_\alpha \sum_{\lambda \in \Y_N^d} \frac{(\dim \V^\lambda)^2}{d!}
			\hat{H}(\alpha)\hat{K}_\alpha(\lambda)
			\hat{\Psi}_\frac{t}{N}(\lambda) \hat{\Omega}_\frac{1}{N}^{-1}(\lambda).
		\end{equation}
		
	\noindent
	The corresponding Gross-Taylor expectation for $1 \leq d \leq N$ is (one hole, one handle)
	
		\begin{equation}
		\label{eqn:HoleyTorusExpectation}
			\mathcal{E}_N^d(\alpha;1) = \langle HK_\alpha \Psi_\frac{t}{N} \Omega_\frac{1}{N}^{-1} \rangle
		\end{equation}
		
		\subsubsection{Microchiral partition function}
		Based on our experience above we can easily predict what will happen in this
		instance of $YM_2$: since spacetime has one hole, the corresponding
		Hurwitz theory must pertain to covers with one free ramification point 
		and simple branching elsewhere; since spacetime has Euler characteristic
		zero the base curve of this Hurwitz theory must be the torus; since the number of holes in spacetime exceeds
		its Euler characteristic this Hurwitz theory must be mixed.
		
		This is borne out by calculations entirely analogous to those already 
		performed. In the stable range, $1 \leq d \leq N$, the expectation 
		\eqref{eqn:HoleyTorusExpectation} admits the absolutely convergent $1/N$ expansion
		
			\begin{equation}
				\mathcal{E}_N^d(\alpha;1) = \sum_{r,s=0}^\infty \frac{(-t)^r (-1)^s}{r!N^{r+s}} \langle HK_\alpha K^r M_s \rangle.
			\end{equation}
		
		\noindent
		Reorganizing this into $\langle HK^rM_sK_\alpha \rangle$, the expectation counts
		tuples $(\rho,\sigma,\tau_1,\dots,\tau_r,\mu_1,\dots,\mu_s,\pi)$ of permutations from 
		$\group{S}(d)$ such that 
		
			\begin{equation}
				\rho^{-1}\sigma^{-1}\rho\sigma\tau_1\dots\tau_r\mu_1\dots\mu_s = \pi.
			\end{equation}
			
		\noindent
		The target permutation $\pi$ belongs to the conjugacy class $K_\alpha$, while 
		$\tau_1,\dots,\tau_r,\mu_1,\dots,\mu_s$ are transpositions with the $\mu$'s weakly 
		monotone increasing. The transpositions $\rho,\sigma$ coming from the commutator 
		sum $|H\rangle$ are arbitrary. Using the Jucys-Murphy correspondence, we can 
		interpret the product $\rho^{-1}\sigma^{-1}\rho\sigma$ as counting a concatenation 
		of four strictly monotone walks, the first two legs of which are strictly monotone 
		decreasing and the final two are strictly monotone increasing, being the first
		two trajectories after path reversal. All together, this gives a walk in the 
		Hurwitz-Cayley graph from $\iota$ to $\pi$ in $2(|\rho|+|\sigma|)+r+s$ steps,
		where the first $2(|\rho|+|\sigma|)$ steps form a walk of shape $\searrow\searrow\nearrow\nearrow$,
		the next $r$ steps are arbitrary, and the final $s$ steps are weakly monotone.
		
		The topological interpretation of the normalized expectation $\frac{1}{d!}\langle HK_\alpha K^r M_s \rangle$
		is that it counts a subset of the covers counted by $\frac{1}{d!}\langle HK_\alpha K^{r+s} \rangle$, 
		which are branched covers of a torus (with no hole) having $1+r+s$ total branch points at fixed
		locations, one of which has ramification type $\alpha$ and the rest are simple.

			\begin{thm}
			\label{thm:MicrochiralOneHoledTorusExpansion}
				For any integers $1 \leq d \leq N$, the microchiral partition function $Z_N^d$
				of $YM_2$ on a one-holed torus of area $t$ admits the absolutely convergent series expansion 
			
					\begin{equation*}
						Z_N^d = \frac{(zN^{-1})^d}{d!} \sum_{r,s=0}^\infty \frac{(-t)^r}{r!N^r} \frac{(-1)^s}{N^s}
						\sum_{\alpha \in \Y^d} p_\alpha \langle HK_\alpha K^r M_s \rangle
					\end{equation*}
				
				\noindent
				where $\langle HK_\alpha K^r M_s \rangle$ is a mixed single Hurwitz number
				with torus target.
			\end{thm}

		\subsubsection{Chiral Gross-Taylor series}
		
			\begin{definition}
				The chiral Gross-Taylor series of the one-hold torus is 
				
					\begin{equation*}
						Z = \sum_{d=0}^\infty \frac{(z\hbar)^d}{d!}\sum_{\alpha \in \Y^d} p_\alpha \sum_{r,s=0}^\infty 
						\frac{(-t\hbar)^r}{r!} (-u\hbar)^s  \langle HK_\alpha K^r M_s \rangle,
					\end{equation*}
					
				\noindent
				where $z,t,u,\hbar,p_1,p_2,p_3,\dots$ are commuting indeterminates.
			\end{definition}
			
			\begin{thm}
				We have $Z= e^F$, where
				
					\begin{equation*}
						F = \sum_{d=0}^\infty \frac{(z\hbar)^d}{d!} \sum_{\alpha \in \Y^d} 
						p_\alpha \sum_{r,s=0}^\infty \frac{ (-t\hbar)^r}{r!} (-u\hbar)^s
						\langle HK_\alpha K^r M_s \rangle
					\end{equation*}
					
				\noindent
				and the cumulant $\langle HK_\alpha K^r M_s \rangle_c$
				is a connected mixed single Hurwitz number with torus target.			
			\end{thm}
			
			The genus expansion is now obtained from the Riemann-Hurwitz formula just
			as in all above computations, e.g. as in the case of a torus without a hole. 
			The quasimodularity of the genus specific free energies $F_g$, $g \geq 1$, which 
			are generating functions for connected mixed single Hurwitz numbers with torus
			target and fixed genus $g$ of the covering surface, is a recent result of \cite{HIL}.
	
	\section{Area Zero}
	\label{sec:AreaZero}
	
		The microchiral partition function of $YM_2$ on an area zero compact orientable surface of with 
		$m$ holes and $n$ handles is 
		
			\begin{equation}
			\label{eqn:AreaZeroGT}
				W_N^d = \frac{\left( zN^{2-2n-m} \right)^d}{d!}
				\sum_{\alpha^1,\dots,\alpha^m \in \Y^d}
				P_{\alpha^1 \dots \alpha^m}
				\sum_{\lambda \in \Y_N^d} \frac{(\dim \V^\lambda)^2}{d!} H^n
				K_{\alpha^1 \dots \alpha^m}(\lambda)\Omega_\frac{1}{N}^{2-2n-m}(\lambda).
			\end{equation}
			
		\noindent
		This is obtained from the true Gross-Taylor formula, Theorem \ref{thm:GrossTaylor}, 
		by setting $t=0$ and it a definition, not a theorem. Not that we have
		retained the variable $z$ in this definition rather than setting it to $1$
		as we would have to were $z=e^{-t/2}$. This is done so that we can define
		the full chiral partition function of $YM_2$ on an area zero surface by 
		
			\begin{equation}
				W_N = 1 + \sum_{d=1}^\infty W_N^d,
			\end{equation}
			
		\noindent
		a formal power series considered in a number of papers in the $YM_2$ literature
		\cite{CMR1,CMR2}. The main result of these works is that in the formal $N=\infty$
		limit, the chiral partition function $W_N$ of area zero $YM_2$ becomes a generating
		function for Euler characteristics of the Hurwitz moduli spaces of branched covers. 
		At the same time, our mixed Hurwitz formalism shows that this
		degeneration of $YM_2$ corresponds exactly to pure monotone Hurwitz theory,
		leading to the conclusion that monotone Hurwitz numbers have an interpretation
		as Euler characteristics of Hurwitz moduli spaces. It would be extremely interesting
		to make this connection precise. Here we give several examples of zero area $YM_2$
		as pure monotone Hurwitz theory.
		
		\subsection{Cylinder}
		The michrochiral partition function of $YM_2$ on a cylinder of area zero is 
		
			\begin{equation}
				W_N^d = N^{2d}\sum_{\alpha,\beta \in \Y^d} P_{\alpha\beta} \sum_{\lambda \in \Y_N^d} \left(\frac{\dim \V^\lambda}{d!}\right)^2
				K_{\alpha\beta}(\lambda).
			\end{equation}
						
		\noindent
		In the stable range, $1 \leq d \leq N$, the internal sum is the Plancherel expectation
		
			\begin{equation}
				\langle K_{\alpha\beta}\rangle = \langle K_\alpha | K_\beta \rangle = \delta_{\alpha\beta} |K_\alpha|.
			\end{equation}
			
		\noindent
		The area zero chiral Gross-Taylor series is therefore
		
			\begin{equation}
				W = 1 + \sum_{d=1}^\infty \frac{z^d}{d!} \sum_{\alpha \in \Y_d} p_\alpha \otimes p_\alpha |K_\alpha|,
			\end{equation}
			
		\noindent
		where the internal sum can be written 
		
			\begin{equation}
				 \sum_{\alpha \in \Y_d} p_\alpha \otimes p_\alpha |K_\alpha| = \sum_{\pi \in \group{S}^d} 
				 p_{t(\pi)} \otimes p_{t(\pi)},
			\end{equation}
			
		\noindent
		where $t(\pi)$ is the cycle type of the permutation $\pi$. This is multiplicative function on permutations,
		and the Exponential Formula gives
		
			\begin{equation}
				W = \exp \sum_{d=1}^\infty \frac{z^d}{d} p_d \otimes p_d.
			\end{equation}
			
		\noindent
		Noting that 
		
			\begin{equation}
				W = 1+\sum_{d=1}^\infty z^d \sum_{\lambda \in \Y^d} s_\lambda \otimes s_\lambda,
			\end{equation} 
			
		\noindent
		we see that this is exactly the Cauchy identity from symmetric function theory.
			
		\subsubsection{Sphere}
		For a sphere of area zero, the chiral Gross-Taylor series is
		
			\begin{equation}
					W= 1 + \sum_{d=1}^\infty \frac{(z\hbar^{-2})^d}{d!}
					\sum_{a,b=0}^{d-1} (u\hbar)^a (v\hbar)^b 
					\langle L_a L_b \rangle.
			\end{equation}
			
		\noindent
		We have that 
		
			\begin{equation}
				\langle L_a L_b \rangle = \langle L_a | L_b \rangle = \delta_{ab} {d \brack d-a},
			\end{equation}
			
		\noindent
		where the right hand side is the Stirling cycle number enumerating permutations
		in $\group{S}^d$ with $d-a$ cycles. This gives 
		
			\begin{equation}
				W= 1 + \sum_{d=1}^\infty \frac{(z\hbar^{-2})^d}{d!}
					\sum_{a=0}^{d-1} \hbar^{2a}  {d \brack d-a}
			\end{equation}
			
		\noindent
		or equivalently 
		
			\begin{equation}
				W = 1 + \sum_{d=1}^\infty \frac{z^d}{d!}
					\sum_{a=1}^d \hbar^{-2a}  {d \brack a} = \exp \left(\hbar^{-2} \sum_{d=1}^\infty \frac{z^d}{d} \right),
			\end{equation}
			
		\noindent
		the standard generating series for Stirling cycle numbers.

		\subsection{Three-holed sphere}
		For a three-holed sphere of area zero, the chiral Gross-Taylor series is 
		
			\begin{equation}
						W= \sum_{d=0}^\infty \frac{(z\hbar)^d}{d!} \sum_{\alpha,\beta,\gamma \in \Y^d} 
						P_{\alpha\beta\gamma} \sum_{s=0}^\infty (-u\hbar)^s
						\langle K_{\alpha\beta\gamma} M_s \rangle,
			\end{equation}
			
		\noindent
		the total generating series for possibly disconnected monotone triple Hurwitz numbers, and
		the corresponding topological expansion is
		
			\begin{equation}
				W = \exp \sum_{g=0}^\infty \hbar^{2g-2} F_g
			\end{equation}
			
		\noindent
		with 
		
			\begin{equation}
					F_g = u^{2g-2}\sum_{d=1}^\infty \frac{z^d}{d!} \sum_{\alpha,\beta,\gamma \in \Y^d}
					P_{\alpha\beta\gamma} (u\hbar)^{\ell(\alpha)+\ell(\beta)+\ell(\gamma)} 
					H_g^{\leq}(\alpha,\beta,\gamma)
				\end{equation}
		
		\noindent
		the generating series for connected monotone triple Hurwitz numbers of genus $g$.
		Interestingly, monotone triple Hurwitz numbers also arise in the context of matrix models, 
		specifically in relation to Jacobi ensembles \cite{GGR}.

\end{document}